# Solution to Morgan's Problem

Qianghui Xiao[a]

[a] *College of Electrical and Information Engineering, Hunan University of Technology, Zhuzhou, Hunan, 412007, PR China*
*E-mail address: 10010@hut.edu.cn;qh.xiao004@163.com*

**Abstract**

In this paper, some preliminaries about Morgan's problem, signal flow graph and controllable linear time-invariant standard system are first introduced in detail. In order to synthesize the necessary and sufficient condition for decoupling system, the first and second necessary conditions, and a sufficient condition for decoupling a controllable linear time-invariant system are secondly analyzed respectively. Therefore, the nonregular static state feedback expression for decoupling system that ensures the internal stability of the uncontrollable subsystem of the decoupled system at the same time is deduced. Then the pole assignment for the controllable subsystem of a decoupled system while ensuring its decoupled state is introduced. Finally, two examples combined with their corresponding signal flow graphs show the simplicity and feasibility of the necessary and sufficient condition for decoupling system described in the paper.

## 1. Introduction

During the last 60 years, a great deal of interest has been brought to the problem of decoupling a linear time-invariant system by static state feedback. Since B. S. Morgan puts forward the row by row decoupling problem known as Morgan's problem (Morgan, 1964; Waren & Mitter, 1975), the first necessary and sufficient condition for the diagonal decoupling problem of square systems was given by Falb and Wolovich (Falb & Wolovich, 1967).

The restrictive case where $m = p+1$ was solved by Glumineau and Moog in the framework of nonlinear systems (Glumineau & Moog, 1992). Morgan's problem without any assumptions on $(F,G)$ was studied by Descusse, Lafay, and Malabre, and they obtained a necessary and sufficient condition for decoupling such systems (Morse, 1973; Descusse, Lafay, & Malabre, 1984, 1988; Descusse, 1987). They believed this result could be applied to general systems. Unfortunately, in general this is not true (Herrera & Lafay, 1993). Herrera and Lafay's research shows two important facts: a) contrary to what is said in (Descusse, Lafay, & Malabre, 1988), a given system and its corresponding shifted system are not equivalent with respect to their static decouplability; and b) the essential orders (Commault, et al., 1986) are not the minimal infinite structure (Commault, Dion, & Torres, 1991; Suda & Umahashi, 1984; Lafay, 2013; Zagalak & Kucera, 2017) which can be reached in a static state feedback decouplable system.

The stability question for decoupled systems is still not be tackled properly for a long time (Wonham & Morse, 1970; Wang, 1992; Eldem, 1996; Wei, Wang, & Cheng, 2010). It remains together with the general solution to the Morgan's problem, an open problem.

The paper is organized as follows: the next section is devoted to some preliminaries in detail. Section 3-5 explain the first and second necessary conditions, and a sufficient condition for decoupling a controllable linear time-invariant system respectively. Section 6 synthesizes the necessary and sufficient condition for decoupling system and then section 7 deduces the nonregular static state feedback expression for decoupling system that ensures the internal stability of the uncontrollable subsystem of a decoupled system at the same time. Finally, section 8 introduces the pole assignment for the controllable subsystem of a decoupled system while ensuring its decoupled state.

In this paper, the nonregular Morgan's problem is considered and analyzed by default, and the analysis of linear time-invariant systems is based on the zero-state response characteristics which will not be mentioned later.

## 2. Preliminaries

Consider a linear time-invariant system $\sum(A,B,C)$ given by

$$\dot{x} = Ax + Bu \quad (1)$$
$$y = Cx \quad (2)$$

where $x \in R^n$, $u \in R^m$, $y \in R^p$, $p \leq m$. $x \triangleq [x_1, x_2, \cdots, x_n]^T$, $u \triangleq [u_1, u_2, \cdots, u_m]^T$, and $y \triangleq [y_1, y_2, \cdots, y_p]^T$. And $A$, $B$, and $C$ are $n \times n$, $n \times m$, and $p \times n$ matrices, respectively. In all what follows, it will be assumed, without any loss of generality, that the system $\sum(A,B,C)$ is right invertible and controllable, $B$ monic and $C$ epic.

The transfer function of the system $\sum(A,B,C)$ is:

$$T(s) = C(sI - A)^{-1}B \quad (3)$$

Morgan's problem can be stated as follows. Under which conditions does there exists a static state feedback:

$$u = Fx + G\upsilon = Fx + G[\upsilon_1, \upsilon_2, \cdots, \upsilon_p]^T \quad (4)$$

such that, for each $i \leq p$, $\upsilon_i$ controls the scalar output $y_i$, without affecting the $(p-1)$ other outputs (Descusse, Lafay, & Malabre, 1988)? In other words, find $(F,G)$ such that the closed-loop transfer function is diagonal (Herrera & Lafay, 1993).

The closed-loop transfer function of the system $\sum(A,B,C)$ is as follows:

$$T_{FG}(s) = C(sI - A - BF)^{-1}BG \quad (5)$$
$$T_{FG}(s) = T(s)(I - F(sI-A)^{-1}B)^{-1}G \quad (6)$$

Without any loss of generality, we shall first consider the case where the closed-loop system is wanted to be integrator decoupled, that is to say: $T_{FG}(s) = diag[s^{-d_1}, s^{-d_2}, \cdots, s^{-d_p}]$, $d_i \in N$. If $m = p$, $G$ is a square matrix and invertible, then equation (4) is called a regular static state feedback and the system is regular, we talk of the regular Morgan's problem. If $m > p$, $G$ is left invertible, then equation (4) is called a nonregular static state feedback and the system is nonregular, we talk of the nonregular Morgan's problem (Descusse, 1991; Zagalak, Lafay, & Herrera, 1993; Herrera & Lafay, 1993).

**Definition 2.1.** The symbol $x$ represents not only the column vector of state variables (denoted as $x = [x_1, x_2, \cdots, x_n]^T$), but also the set of state variables (denoted as $x = \{x_i\}_n$), and at the same time the set of first-order branches in the signal flow graph.

**Definition 2.2.** The symbol $x_i$ represents not only the state variable, but also the first-order branch with $x_i$ as its output, and at the same time the output of the first-order branch in the signal flow graph.

**Definition 2.3.** The symbol $y$ represents not only the column vector of outputs (denoted as $y = [y_1, y_2, \cdots, y_p]^T$), but also the set of outputs (denoted as $y = \{y_i\}_p$).

**Definition 2.4.** The symbol $u$ represents not only the column vector of inputs (denoted as $u = [u_1, u_2, \cdots, u_m]^T$), but also the set of control inputs (denoted as $u = \{u_i\}_m$).

### 2.1. Some basic concepts about Morgan's problem

**Definition 2.5.** Let $d_1, d_2, \cdots, d_p$ be given by

$$d_i = \min\{j : C_i A^{j-1} B \neq 0, j = 1, 2, \cdots, n\} \quad (7)$$

or

$$d_i = n \quad \text{if } C_i A^{j-1} B = 0 \text{ for all } j \quad (8)$$

$$C \triangleq [C_1^T, C_2^T, \cdots, C_p^T]^T \quad (9)$$

where $C_i$ denotes the *ith* row of $C$, and integer $d_i$ is defined as the relative order corresponding to the output $y_i$.

**Definition 2.6.** Let $B^*$ be the $p \times m$ matrix given by

$$B^* = \begin{bmatrix} C_1 A^{d_1-1} B \\ C_2 A^{d_2-1} B \\ \vdots \\ C_p A^{d_p-1} B \end{bmatrix} \triangleq \begin{bmatrix} B_1^* \\ B_2^* \\ \vdots \\ B_p^* \end{bmatrix} \quad (10)$$

$$B^* \triangleq [\bar{B}_1^* \quad \bar{B}_2^* \quad \cdots \quad \bar{B}_m^*] \quad (11)$$

Let $A^*$ be the $p \times n$ matrix given by

$$A^* = \begin{bmatrix} C_1 A^{d_1} \\ C_2 A^{d_2} \\ \vdots \\ C_p A^{d_p} \end{bmatrix} \triangleq \begin{bmatrix} A_1^* \\ A_2^* \\ \vdots \\ A_p^* \end{bmatrix} \quad (12)$$

**Proposition 2.7.** For a regular Morgan's problem, there exists a static state feedback $(F,G)$ which decouple the system $\sum(A,B,C)$ if and only if

$$\det B^* \neq 0 \quad (13)$$

i.e., if and only if $B^*$ is invertible or nonsingular (Falb & Wolovich, 1967).

$$u = Fx + G\upsilon \triangleq Fx + G[\upsilon_1, \upsilon_2, \cdots, \upsilon_p]^T \quad (14)$$
$$\upsilon = A^* x + B^* u \quad (15)$$
$$F = -(B^*)^{-1} A^* \quad (16)$$
$$G = (B^*)^{-1} \quad (17)$$
$$\dot{x} = (A + BF)x + BG\upsilon \triangleq \hat{A}x + \hat{B}\upsilon \quad (18)$$

$\upsilon_1, \upsilon_2, \cdots, \upsilon_p$ are defined as the feedforward (control) inputs in the static state feedback $(F,G)$.

It is assumed that the set of the remaining first-order branches except the controllable and observable first-order branches of the decoupled system $\sum(\hat{A}, \hat{B}, C)$ is $x'$, then an equivalent decoupled (tree-like) system which is used to analyze the pole assignment of the system $\sum(\hat{A}, \hat{B}, C)$ can be derived as follows.

$$\hat{x}^T \triangleq \left[ Y^T, (x')^T \right] \quad (19)$$
$$Y^T \triangleq \left[ Y_1^T, Y_2^T, \cdots, Y_p^T \right] \quad (20)$$
$$Y_i^T \triangleq \left[ y_i, y_i^{(1)}, \cdots, y_i^{(d_i-1)} \right], i = 1, 2, \cdots, p \quad (21)$$
$$Y_i = P_i x, i = 1, 2, \cdots, p \quad (22)$$
$$P_i \triangleq \begin{bmatrix} C_i \\ C_i \hat{A} \\ \vdots \\ C_i \hat{A}^{d_i-1} \end{bmatrix}, i = 1, 2, \cdots, p \quad (23)$$
$$P \triangleq \left[ P_1^T, P_2^T, \cdots, P_p^T \right]^T \triangleq \left[ \bar{P}_1 \quad \bar{P}_2 \right] \quad (24)$$

$$n_{co} \triangleq \sum_{i=1}^{p} d_i \tag{25}$$

Since the $n_{co}$ variables in $Y$ are linearly independent, the matrix $P$ is row-full-rank.

$$rankP = n_{co} \tag{26}$$

$$\hat{x} = Tx \triangleq \begin{bmatrix} \bar{P}_1 & \bar{P}_2 \\ 0 & I \end{bmatrix} x \tag{27}$$

$$\det T \neq 0 \tag{28}$$

$$\begin{cases} \dot{\hat{x}} = T\hat{A}T^{-1}\hat{x} + T\hat{B}\upsilon \triangleq \tilde{A}\hat{x} + \tilde{B}\upsilon \\ y = CT^{-1}\hat{x} \triangleq \tilde{C}\hat{x} \end{cases} \tag{29}$$

$$y_i^{(d_i)} = \upsilon_i, i = 1, 2, \cdots, p \tag{30}$$

$$\dot{Y}_i = \tilde{A}_{ii} Y_i + b_i \upsilon_i, i = 1, 2, \cdots, p \tag{31}$$

$$\tilde{A}_{ii} = \begin{bmatrix} 0 & 1 & 0 & \cdots & 0 \\ 0 & \ddots & \ddots & \ddots & \vdots \\ \vdots & \ddots & \ddots & \ddots & 0 \\ 0 & \cdots & 0 & 0 & 1 \\ 0 & 0 & 0 & 0 & 0 \end{bmatrix} \tag{32}$$

$$b_i = \begin{bmatrix} 0 & 0 & \cdots & 0 & 1 \end{bmatrix}^T \tag{33}$$

$$T_{FG}(s) = diag[s^{-d_1}, s^{-d_2}, \cdots, s^{-d_p}] \tag{34}$$

The equations (31) represent that the controllable and observable subsystem of the equivalent decoupled system $\sum(\tilde{A}, \tilde{B}, \tilde{C})$ is only composed of $p$ independent chains of integrators with each order of $d_i$ (i.e., $p$ mutually independent forward paths). But it cannot be guaranteed that the controllable subsystem of $\sum(\tilde{A}, \tilde{B}, \tilde{C})$ must be composed of $p$ mutually independent controllable subsystems.

**Theorem 2.8.** For a regular Morgan's problem, let

$$\hat{\upsilon} = B^* u \tag{35}$$

$$\hat{\upsilon} \triangleq \{\hat{\upsilon}_i\}_p \triangleq [\hat{\upsilon}_1, \hat{\upsilon}_2, \cdots, \hat{\upsilon}_p]^T \tag{36}$$

Then there exists a static state feedback $(F, G)$ which decouple the system $\sum(A, B, C)$ if and only if $\hat{\upsilon}_1, \hat{\upsilon}_2, \cdots, \hat{\upsilon}_p$ are linearly independent that is equivalent to $\det B^* \neq 0$.

**Proof.** $u_1, u_2, \cdots, u_p$ are mutually independent control inputs, and undetermined coefficient is assumed to be $\alpha \in C^{p \times 1}$, $\alpha \triangleq [\alpha_1, \alpha_2, \cdots, \alpha_p]^T$,
Regular Morgan's problem is solvable
$\Leftrightarrow \det B^* \neq 0 \Leftrightarrow \alpha^T B^* = 0$, if and only if $\alpha = 0$
$\Leftrightarrow \alpha^T B^* u = 0$, if and only if $\alpha = 0$
$\Leftrightarrow \alpha^T \hat{\upsilon} = 0$, if and only if $\alpha = 0$
$\Leftrightarrow \hat{\upsilon}_1, \hat{\upsilon}_2, \cdots, \hat{\upsilon}_p$ are linearly independent
□

**Theorem 2.9.** For a nonregular Morgan's problem, let

$$\hat{\upsilon} = B^* u \tag{37}$$

$$\hat{\upsilon} \triangleq \{\hat{\upsilon}_i\}_p \triangleq [\hat{\upsilon}_1, \hat{\upsilon}_2, \cdots, \hat{\upsilon}_p]^T \tag{38}$$

Then there exists a static state feedback $(F, G)$ which decouple the system $\sum(A, B, C)$ if $\hat{\upsilon}_1, \hat{\upsilon}_2, \cdots, \hat{\upsilon}_p$ are linearly independent that is equivalent to $rankB^* = p$.

**Proof.** $rankB^* = p$, there is no problem to assume that the first $p$ columns of the matrix $B^*$ are linearly independent.

$$\bar{B}^* \triangleq [\bar{B}_1^* \quad \bar{B}_2^* \quad \cdots \quad \bar{B}_p^*] \tag{39}$$

$$\tilde{B}^* \triangleq [\bar{B}_{p+1}^* \quad \bar{B}_{p+2}^* \quad \cdots \quad \bar{B}_m^*] \tag{40}$$

$$U_1 \triangleq [u_1, u_2, \cdots, u_p]^T \tag{41}$$

$$U_2 \triangleq [u_{p+1}, u_{p+2}, \cdots, u_m]^T \tag{42}$$

$$\det \bar{B}^* \neq 0 \tag{43}$$

$$\bar{B}_j^* = \sum_{i=1}^{p} \beta_{ij} \bar{B}_i^*, p < j \leq m \tag{44}$$

$$\upsilon = A^* x + B^* u = A^* x + \bar{B}^* U_1 + \tilde{B}^* U_2 \tag{45}$$

$$U_2 = K_2 x \tag{46}$$

$$\upsilon = (A^* + \tilde{B}^* K_2) x + \bar{B}^* U_1 \tag{47}$$

$$U_1 = -(\bar{B}^*)^{-1}(A^* + \tilde{B}^* K_2) x + (\bar{B}^*)^{-1} \upsilon \tag{48}$$

$$F = \begin{bmatrix} -(\bar{B}^*)^{-1}(A^* + \tilde{B}^* K_2) \\ K_2 \end{bmatrix} \tag{49}$$

$$G = \begin{bmatrix} (\bar{B}^*)^{-1} \\ 0 \end{bmatrix} \tag{50}$$

$$u = Fx + G\upsilon \tag{51}$$

$$T_{FG}(s) = diag[s^{-d_1}, s^{-d_2}, \cdots, s^{-d_p}] \tag{52}$$

$u_1, u_2, \cdots, u_p$ are called master (control) inputs of the system $\sum(A, B, C)$, and $u_{p+1}, u_{p+2}, \cdots, u_m$ are called auxiliary (control) inputs. (46) is the self-compensation of all auxiliary inputs of $\sum(A, B, C)$ explained in detail later.

$u_1, u_2, \cdots, u_p$ are mutually independent controls, and undetermined coefficient is assumed to be $\alpha \in C^{p \times 1}$, $\alpha \triangleq [\alpha_1, \alpha_2, \cdots, \alpha_p]^T$,
$rankB^* = p \Leftrightarrow \det \bar{B}^* \neq 0$
$\Leftrightarrow \alpha^T \bar{B}^* = 0$, if and only if $\alpha = 0$
$\Leftrightarrow \alpha^T B^* = 0$, if and only if $\alpha = 0$
$\Leftrightarrow \alpha^T B^* u = 0$, if and only if $\alpha = 0$
$\Leftrightarrow \alpha^T \hat{\upsilon} = 0$, if and only if $\alpha = 0$
$\Leftrightarrow \hat{\upsilon}_1, \hat{\upsilon}_2, \cdots, \hat{\upsilon}_p$ are linearly independent
□

**Proposition 2.10.** Coordinate transformation does not change the controllability, observability and row by row decouplability of the system $\sum(A, B, C)$ (Wolovich, 1974; Wonham, 1974; Pang, 1992).

**Proposition 2.11.** Static state feedback does not change the controllability and row by row decouplability of $\sum(A, B, C)$ (Wolovich, 1974; Wonham, 1974; Pang, 1992).

**Proposition 2.12.** For a nonregular Morgan's problem, right invertibility of $\sum(A, B, C)$ is a necessary condition for its row by row decouplability, that is, $T(s)$ is right invertible with the rank of $p$ (Descusse, Lafay, & Malabre, 1988).

**Proposition 2.13.** There exists a general solution of Morgan's problem for $\sum(A, B, C)$ if and only if there exists a control law $(F_0, G_0)$ such that $\sum(A + BF_0, BG_0, C)$ is decoupleable (Descusse, Lafay, & Malabre, 1988).

*2.2. Fundamentals about signal flow graph*

The signal flow graph originated from Mason's use of graphical method to describe linear algebraic equations. It is a signal transmission network composed of nodes and directed branches (referred as branches), and is a typical application of directed graph (or digraph). Based on the zero-state response, a linear time-invariant system can be converted into a corresponding linear algebraic system by using Laplace transform method, so that a signal flow graph can be used to describe the linear time-invariant system. Linear time-invariant system and signal flow graph can be converted to each other, so that there is a one-to-one correspondence between the two (Bang & Gutin, 2002).

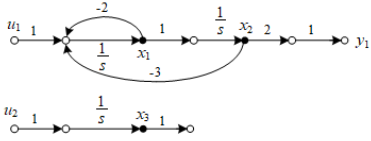

Fig.1. A schematic signal flow graph

It is assumed that the branches in signal flow graphs are divided into two types: general first-order branch and zero-order branch with a constant gain. One general first-order branch consists of a first-order branch (its input node is the input node of the general first-order branch, and its output is the state variable corresponding to this first-order branch) with a gain of $1/s$ and a zero-order branch with a constant gain (i.e., combination coefficient of this first-order branch) in series. A schematic signal flow graph is shown in Fig.1.

**Definition 2.14.** A route from a node (i.e., input node) such as $z_1$ to another node (i.e., output node) such as $z_2$ is called a walk (denoted as $walk(z_1, z_2)$), and a walk that passes through its each node only once is called a path (denoted as $path(z_1, z_2)$). One path from a control input node such as $u_i$ to an output node such as $y_j$ is called a forward path (denoted as $path(u_i, y_j)$). Multiple paths without common nodes between each other are called node-disjoint paths or non-contact paths. A path is called a chain of integrators (or integrator-string) while its input node is a control input and its output node is a certain state variable.

**Definition 2.15.** A closed path with its start and end points at the same node is called a loop or cycle. Multiple loops without common nodes between each other are called node-disjoint loops or non-contact loops. A loop that intersects a first-order branch with its input node is called a contact loop or node-joint loop of this first-order branch. Then a loop that intersects a first-order branch with its output node is called a branch-joint loop of this first-order branch. A loop set is called an associative loop set when each first-order branch in this loop set is connected to another first-order branch by one path in the loop set. An associative loop set is called the associative loop set of a certain first-order branch when the first-order branch and this associative loop set are branch-joint. If no node intersects between two loop sets (especially two associative loop sets), they are called two node-disjoint loop sets (especially two node-disjoint associative loop sets).

**Definition 2.16.** The order or length of a walk (or path, or loop, or chain of integrators) is defined as the number of first-order branches contained in this walk (or path, or loop, or chain of integrators).

**Proposition 2.17.** For a linear time-invariant system, if there is at least one forward path from a control input node (such as $u_i$) to an output node (such as $y_j$), then there exists at least one path with the shortest length (i.e., shortest path, denoted as $path(u_i, y_j)_{\min}$).

**Definition 2.18.** A signal flow graph without loops is called a tree-like signal flow graph (corresponding to its tree-like linear system).

**Proposition 2.19.** The open-loop and closed-loop poles of each first-order branch (i.e., its corresponding state variable) in a tree-like signal flow graph are zero.

**Proposition 2.20.** After a static state feedback for a linear time-invariant system, as long as a first-order branch in the corresponding signal flow graph increases or decreases at least one loop in its associative loop set (that is, its associative loop set has changed), the closed-loop pole of this first-order branch will change. This means that the closed-loop pole of each first-order branch is determined by its associative loop set.

It can be seen that a signal flow graph with the open-loop pole of each first-order branch at zero can well show how static state feedback changes the closed-loop pole of each first-order branch.

In the corresponding signal flow graph for a linear time-invariant system, each state variable (their column vector or set such as $X_1$) corresponding to the first-order branch that there is at least one path connected from a control input (such as $u_i$) to it can construct the subsystem controlled by this control input. The column vector or set $X_1^{'}$ is defined as: $X_1 \cap X_1^{'} = \varnothing$ and $X_1 \cup X_1^{'} = x$. According to the linear superposition principle, other control inputs except $u_i$ are set to zero, and $X_1^{'} = 0$ can be obtained based on the zero state response to construct the subsystem controlled by $u_i$.

$$\begin{bmatrix} \dot{X}_1 \\ \dot{X}_1^{'} \end{bmatrix} = \begin{bmatrix} A_{11} & A_{12} \\ 0 & A_{22} \end{bmatrix} \begin{bmatrix} X_1 \\ X_1^{'} \end{bmatrix} + \begin{bmatrix} B_1 & B_2 & \cdots & B_p \end{bmatrix} u \quad (53)$$

$$X_1^{'} = 0 \quad (54)$$

$$u_j = 0, j \neq i \quad (55)$$

$$\dot{X}_1 = A_{11} X_1 + B_{i1} u_i \quad (56)$$

$$B_i^T \triangleq [B_{i1}^T, B_{i2}^T] \quad (57)$$

$$y = Cx \triangleq \begin{bmatrix} C_1 & C_1^{'} \end{bmatrix} \begin{bmatrix} X_1 \\ X_1^{'} \end{bmatrix} = C_1 X_1 \quad (58)$$

**Definition 2.21.** For a controllable linear time-invariant system $\Sigma(A, B, C)$, the controllable subsystem controlled by the control input $u_i$ is $\Sigma(A_{11}, B_{i1})$.

Correspondingly, in the case of multiple control inputs (such as $U_1 \triangleq [u_{i_1}, u_{i_2}, \cdots, u_{i_p}]^T$, and the column vector $U_2$ formed by the remaining control inputs), the column vector or set of the state variables (the elements of which are different from each other) obtained by combining the state variables controlled by each control input is $\bar{X}_1$, and the column vector or set $X_2$ is defined as: $\bar{X}_1 \cap X_2 = \varnothing$ and $\bar{X}_1 \cup X_2 = x$. The controllable subsystem controlled by the multiple control inputs can be obtained for the same reason as above.

$$\begin{bmatrix} \dot{\bar{X}}_1 \\ \dot{X}_2 \end{bmatrix} = \begin{bmatrix} \bar{A}_{11} & \bar{A}_{12} \\ 0 & \bar{A}_{22} \end{bmatrix} \begin{bmatrix} \bar{X}_1 \\ X_2 \end{bmatrix} + \begin{bmatrix} \bar{B}_1 & \bar{B}_2 \end{bmatrix} \begin{bmatrix} U_1 \\ U_2 \end{bmatrix} \quad (59)$$

$$X_2 = 0 \tag{60}$$

$$U_2 = 0 \tag{61}$$

$$\dot{\bar{X}}_1 = \bar{A}_{11}\bar{X}_1 + \bar{B}_{11}U_1 \tag{62}$$

$$\bar{B}_1^T \triangleq \begin{bmatrix} \bar{B}_{11}^T & \bar{B}_{12}^T \end{bmatrix} \tag{63}$$

$$y = Cx \triangleq \begin{bmatrix} C_1 & C_2 \end{bmatrix} \begin{bmatrix} \bar{X}_1 \\ X_2 \end{bmatrix} = C_1 \bar{X}_1 \tag{64}$$

**Definition 2.22.** For a controllable linear time-invariant system $\Sigma(A,B,C)$, the controllable subsystem controlled by multiple control inputs (such as $U_1$) is $\Sigma(\bar{A}_{11}, \bar{B}_{11})$, and the uncontrollable part is $X_2$ corresponding to $U_1$.

Firstly, we can only assign the closed-loop poles for the uncontrollable part $X_2$ corresponding to $U_1$ by the static state feedback $(F_{21}, I)$ as follows.

$$U_2 = F_{21}X_2 + \bar{U}_2 \tag{65}$$

$$U_2 = \begin{bmatrix} 0 & F_{21} \end{bmatrix} \begin{bmatrix} \bar{X}_1 \\ X_2 \end{bmatrix} + \bar{U}_2 \tag{66}$$

$$U_2 = F_{22}x + \bar{U}_2 \tag{67}$$

$$u = \begin{bmatrix} U_1 \\ U_2 \end{bmatrix} = \begin{bmatrix} 0 \\ F_{22} \end{bmatrix} x + \begin{bmatrix} U_1 \\ \bar{U}_2 \end{bmatrix} \triangleq F_1 x + \begin{bmatrix} U_1 \\ \bar{U}_2 \end{bmatrix} \tag{68}$$

Secondly, a static state feedback $(F_{11}, I)$ can only be designed to assign each closed-loop pole for the controllable part $\bar{X}_1$ (including the controllable and unobservable part) corresponding to $U_1$ as above to have negative real part, so that the controllable and unobservable part can be transformed into the controllable and observable part.

$$U_1 = F_{11}\bar{X}_1 + \bar{U}_1 \tag{69}$$

$$U_1 = \begin{bmatrix} F_{11} & 0 \end{bmatrix} \begin{bmatrix} \bar{X}_1 \\ X_2 \end{bmatrix} + \bar{U}_1 \tag{70}$$

$$U_1 \triangleq F_{12}x + \bar{U}_1 \tag{71}$$

$$u = \begin{bmatrix} U_1 \\ U_2 \end{bmatrix} = \begin{bmatrix} F_{12} \\ 0 \end{bmatrix} x + \begin{bmatrix} \bar{U}_1 \\ U_2 \end{bmatrix} = \begin{bmatrix} F_{12} \\ F_{22} \end{bmatrix} x + \begin{bmatrix} \bar{U}_1 \\ \bar{U}_2 \end{bmatrix} \tag{72}$$

$$\bar{u} \triangleq \begin{bmatrix} \bar{U}_1 \\ \bar{U}_2 \end{bmatrix} \tag{73}$$

$$u = F_2 x + \bar{u} \triangleq \begin{bmatrix} F_{11} & 0 \\ 0 & F_{21} \end{bmatrix} \begin{bmatrix} \bar{X}_1 \\ X_2 \end{bmatrix} + \bar{u} \tag{74}$$

where $\bar{U}_1$ (or $\bar{U}_2$) is the feedforward (control) input corresponding to $U_1$ (or $U_2$).

**Theorem 2.23.** For a controllable linear time-invariant system $\Sigma(A,B,C)$, the static state feedbacks $(F_{21}, I)$ and $(F_{11}, I)$ are equivalent to a static state feedback $(F_2, I)$ for system $\Sigma(A,B,C)$.

**Definition 2.24.** The self-compensation of $U_2$ is expressed as (75) by setting $\hat{U}_2 = 0$ in equation (67). It is equivalent to inserting a chain of integrators with infinite length (or order) before each control input in $U_2$.

$$U_2 = F_{22}x \tag{75}$$

$$\dot{x} = Ax + \begin{bmatrix} \bar{B}_1 & \bar{B}_2 \end{bmatrix} \begin{bmatrix} U_1 \\ U_2 \end{bmatrix} = Ax + \bar{B}_1 U_1 + \bar{B}_2 U_2 \tag{76}$$

$$\dot{x} = (A + \bar{B}_1 F_{12} + \bar{B}_2 F_{22})x + \bar{B}_1 \bar{U}_1 \triangleq \bar{A}_1 x + \bar{B}_1 \bar{U}_1 \tag{77}$$

If $U_1 \in C^{p \times 1}$, after the self-compensation of $U_2$, the system $\Sigma(A,B,C)$ is transformed into a regular system $\Sigma(\bar{A}_1, \bar{B}_1, C)$ and then the uncontrollable subsystem of $\Sigma(\bar{A}_1, \bar{B}_1, C)$ may be produced.

**Definition 2.25.** $k$-path problem can be stated as: for a given signal flow graph and its different nodes such as $u_1, u_2, \cdots, u_k$, $y_1, y_2, \cdots, y_k$, whether there is a set of node-disjoint paths such as $P_1, P_2, \cdots, P_k$ in this signal flow graph, so that each path $P_i$ is $path(u_i, y_i), i = 1, 2, \cdots, k$ (Bang & Gutin, 2002)?

**Definition 2.26.** General $k$-path problem can be stated as: for a given signal flow graph and its different nodes such as $u_1, u_2, \cdots, u_k$, $y_1, y_2, \cdots, y_k$, whether there is a set of node-disjoint paths such as $P_1, P_2, \cdots, P_k$ (called a general $k$-path) in this signal flow graph, so that each path $P_i$ (called a general path $P_i$) is $path(u_{j_i}, y_i)$, $i = 1, 2, \cdots, k$? The sequence $j_1, j_2, \cdots, j_k$ is a kind of mathematic permutation of $1, 2, \cdots, k$.

**Definition 2.27.** If $r > k$, general $k$-path problem can also be stated as: for a given signal flow graph and its different nodes such as $u_1, u_2, \cdots, u_r$, $y_1, y_2, \cdots, y_k$, whether there is a set of node-disjoint paths such as $P_1, P_2, \cdots, P_k$ (called a general $k$-path) in the signal flow graph, so that each path $P_i$ (called a general path $P_i$) is $path(u_{j_i}, y_i)$, $i = 1, 2, \cdots, k$? The sequence $j_1, j_2, \cdots, j_r$ is a kind of mathematic permutation of $1, 2, \cdots, r$.

The general $k$-path problem is a generalized treatment for a $k$-path problem with the different mathematic permutations for its control input indexes.

### 2.3. Controllable linear time-invariant standard system

**Proposition 2.28.** Assuming that the eigenvalues of a $n \times n$ matrix $A$ are $p_1, p_2, \cdots, p_n$ and its characteristic polynomial is $\Delta(s)$, then the following conclusions hold.

$$\Delta(s) = \prod_{i=1}^{n}(s - p_i) \triangleq \sum_{i=0}^{n} \alpha_i s^i, \alpha_n = 1 \tag{78}$$

$$\Delta(A) = \sum_{i=0}^{n} \alpha_i A^i = \prod_{i=1}^{n}(A - p_i I) = 0 \tag{79}$$

**Proposition 2.29.** For a SISO controllable linear time-invariant system $\Sigma(A,b,c)$, there exists one and only one static state feedback as shown in (81) to assign all closed-loop poles of the system as follows (Wolovich, 1974; Wonham, 1974; Wonham, 1979; Pang, 1992; Dion, Commault, &Woude, 2003).

$$\begin{cases} \dot{x} = Ax + bu \\ y = cx \end{cases} \tag{80}$$

$$u = Kx + \bar{u} \tag{81}$$

$$\dot{x} = (A + bK)x + b\bar{u} \triangleq \bar{A}x + b\bar{u} \tag{82}$$

It is assumed that the designated closed-loop poles of the system are $p_1, p_2, \cdots, p_n$ and the characteristic polynomial of $n \times n$ matrix $\bar{A}$ is $\bar{\Delta}(s)$, the following conclusions can be obtained by further derivation.

$$\bar{\Delta}(s) = \prod_{j=1}^{n}(s - p_j) \triangleq \sum_{j=0}^{n} \alpha_j s^j, \alpha_n = 1 \tag{83}$$

$$\bar{\Delta}(\bar{A}) = \sum_{j=0}^{n} \alpha_j \bar{A}^j = \prod_{j=1}^{n}(\bar{A} - p_j I) = 0 \quad (84)$$

$$y^{(j)} = c\bar{A}^j x, 0 \le j < d \quad (85)$$

$$y^{(d)} = c\bar{A}^d x + c\bar{A}^{d-1} b\bar{u} \quad (86)$$

$$y^{(j)} = c\bar{A}^j x + \sum_{k=d}^{j} c\bar{A}^{k-1} b\bar{u}^{(k-d)}, d \le j \le n \quad (87)$$

$$d = \min\{k : cA^{k-1}b \ne 0, k = 1, 2, \cdots, n\} \quad (88)$$

or

$$d = n \quad \text{if } cA^k b = 0 \text{ for all } k \quad (89)$$

where integer $d$ is defined as the relative order corresponding to the output $y$.

$$\sum_{j=0}^{n} \alpha_j y^{(j)} = \sum_{j=0}^{n} \alpha_j c\bar{A}^j x + \sum_{j=d}^{n} (\sum_{k=d}^{j} \alpha_k c\bar{A}^{k-1} b\bar{u}^{(k-d)}) \quad (90)$$

$$\sum_{j=0}^{n} \alpha_j s^j y(s) = c(\sum_{j=0}^{n} \alpha_j \bar{A}^j)x(s) + \sum_{j=d}^{n}\sum_{k=d}^{j} \alpha_k c\bar{A}^{k-1} bs^{k-d}\bar{u}(s) \quad (91)$$

$$\bar{\Delta}(s) y(s) = c\bar{\Delta}(\bar{A})x(s) + \sum_{j=d}^{n}\sum_{k=d}^{j} \alpha_k c\bar{A}^{k-1} bs^{k-d}\bar{u}(s) \quad (92)$$

$$m(s) \triangleq \sum_{j=d}^{n}\sum_{k=d}^{j} \alpha_k c\bar{A}^{k-1} bs^{k-d} \quad (93)$$

$$\bar{\Delta}(s) y(s) = c\bar{\Delta}(\bar{A})x(s) + m(s)\bar{u}(s) = m(s)\bar{u}(s) \quad (94)$$

$$T_F \triangleq y(s)/\bar{u}(s) = m(s)/\bar{\Delta}(s) \quad (95)$$

**Definition 2.30.** A tree-like controllable linear time-invariant system (corresponding to its tree-like signal flow graph) is defined as the controllable linear system whose all open-loop and closed-loop poles are zero, which can also be called a controllable linear time-invariant standard system specially for decoupling problem.

**Proposition 2.31.** For a controllable linear time-invariant system $\sum(A, B, C)$, then an equivalent transformation matrix $T_{c2}$ can be obtained and $\bar{x} = T_{c2} x$ transforms the system into a Luenberger second controllable canonical form as follows (Wolovich, 1974; Wonham, 1974; Pang, 1992).

$$\dot{\bar{x}} = A_{c2}\bar{x} + B_{c2}u \quad (96)$$

$$y = C_{c2}\bar{x} \quad (97)$$

where

$$A_{c2} = \begin{bmatrix} A_{11} & A_{12} & \cdots & A_{1m} \\ A_{21} & A_{22} & \cdots & A_{2m} \\ \vdots & \vdots & \ddots & \vdots \\ A_{m1} & A_{m2} & \cdots & A_{mm} \end{bmatrix} \triangleq \begin{bmatrix} A_1 \\ A_2 \\ \vdots \\ A_m \end{bmatrix}, A_{ij} \in C^{\sigma_i \times \sigma_j} \quad (98)$$

$$A_{ii} = \begin{bmatrix} 0 & 1 & 0 & \cdots & 0 \\ 0 & \ddots & \ddots & \ddots & \vdots \\ \vdots & \ddots & \ddots & \ddots & 0 \\ 0 & \cdots & 0 & 0 & 1 \\ * & * & * & * & * \end{bmatrix} \quad (99)$$

$$A_{ij} = \begin{bmatrix} 0 & 0 & \cdots & 0 \\ \vdots & \vdots & \ddots & \vdots \\ 0 & 0 & \cdots & 0 \\ * & * & \cdots & * \end{bmatrix}, j \ne i \quad (100)$$

$$B_{c2}^T = \begin{bmatrix} B_1^T & B_2^T & \cdots & B_m^T \end{bmatrix}, B_i \in C^{\sigma_i \times m} \quad (101)$$

$$B_i(k, j) = \begin{cases} 0, j < i, k = \sigma_i \\ 1, j = i, k = \sigma_i \\ \beta_{ij}, j > i, k = \sigma_i \\ 0, \forall j, k \ne \sigma_i \end{cases}, i, j = 1, 2, \cdots, m \quad (102)$$

$$C_{c2} = CT_{c2}^{-1} \quad (103)$$

The symbol "*" in above matrices represents the entry that is possibly nonzero and $\beta_{ij}(j > i)$ are undetermined coefficients. $\sigma_1, \sigma_2, \cdots, \sigma_m$ are called the controllable structure indexes of the system, and $\sigma_1 + \sigma_2 + \cdots + \sigma_m = n$ for a controllable linear time-invariant system.

In the following, a regular static state feedback as shown in (112) is used to further transform the above system $\sum(A_{c2}, B_{c2}, C_{c2})$ into a better applicable controllable standard system with no loop structure.

Set

$$\bar{u} \triangleq \begin{bmatrix} \bar{u}_1 & \bar{u}_2 & \cdots & \bar{u}_m \end{bmatrix}^T \quad (104)$$

$$\bar{u}_m = u_m + A_{m;\sigma_m}\bar{x} \quad (105)$$

$$\bar{u}_i = u_i + \sum_{j=i+1}^{m} \beta_{ij} u_j + A_{i;\sigma_i}\bar{x}, i = m, m-1, \cdots, 1 \quad (106)$$

$$\bar{u} = A^*\bar{x} + B^* u \quad (107)$$

$$A^* = \begin{bmatrix} A_{1;\sigma_1} \\ A_{2;\sigma_2} \\ \vdots \\ A_{m;\sigma_m} \end{bmatrix}, A_{i;\sigma_i} \in C^{1 \times n}, i = 1, 2, \cdots, m \quad (108)$$

$$B^* = \begin{bmatrix} B_{1;\sigma_1} \\ B_{2;\sigma_2} \\ \vdots \\ B_{m;\sigma_m} \end{bmatrix}, B_{i;\sigma_i} \in C^{1 \times m}, i = 1, 2, \cdots, m \quad (109)$$

$$B^* = \begin{bmatrix} B_{1;\sigma_1} \\ B_{2;\sigma_2} \\ \vdots \\ B_{m;\sigma_m} \end{bmatrix} = \begin{bmatrix} 1 & \beta_{12} & \beta_{13} & \cdots & \beta_{1m} \\ 0 & 1 & \beta_{23} & \cdots & \beta_{2m} \\ \vdots & 0 & 1 & \ddots & \vdots \\ 0 & 0 & \ddots & \ddots & \beta_{(m-1)m} \\ 0 & 0 & 0 & 0 & 1 \end{bmatrix} \quad (110)$$

$$\det B^* = 1 \ne 0 \quad (111)$$

$$u = F_{c3}\bar{x} + G_{c3}\bar{u} = F_{c3}T_{c2}x + G_{c3}\bar{u} \triangleq K_{c3}x + G_{c3}\bar{u} \quad (112)$$

$$F_{c3} = -(B^*)^{-1} A^* \quad (113)$$

$$G_{c3} = (B^*)^{-1} \quad (114)$$

where $A_{i;\sigma_i}$ (or $B_{i;\sigma_i}$) is the $\sigma_i$th row of $\sigma_i \times n$ matrix $A_i$ (or $\sigma_i \times m$ matrix $B_i$).

So

$$\overline{x}_{\overline{\sigma}_i}^{(\sigma_i)} = \overline{u}_i, i=1,2,\cdots,m \qquad (115)$$

$$\overline{\sigma}_i \triangleq \sum_{j=1}^{i-1} \sigma_j + 1, i=1,2,\cdots,m \qquad (116)$$

$$\dot{\overline{x}} = A_{c3}\overline{x} + B_{c3}\overline{u} \qquad (117)$$

$$A_{c3} = diag\begin{bmatrix} \overline{A}_{11} & \overline{A}_{22} & \cdots & \overline{A}_{mm} \end{bmatrix}, \overline{A}_{ii} \in C^{\sigma_i \times \sigma_i} \qquad (118)$$

$$\overline{A}_{ii} = \begin{bmatrix} 0 & 1 & 0 & \cdots & 0 \\ 0 & \ddots & \ddots & \ddots & \vdots \\ \vdots & \ddots & \ddots & \ddots & 0 \\ 0 & \cdots & 0 & 0 & 1 \\ 0 & 0 & 0 & 0 & 0 \end{bmatrix} \qquad (119)$$

$$B_{c3}^T = \begin{bmatrix} \overline{B}_1^T & \overline{B}_2^T & \cdots & \overline{B}_m^T \end{bmatrix}, \overline{B}_i \in C^{\sigma_i \times m} \qquad (120)$$

$$\overline{B}_i(k,j) = \begin{cases} 1, j=i, k=\sigma_i \\ 0, otherwise \end{cases}, i,j=1,2,\cdots,m \qquad (121)$$

$$y = C_{c2}\overline{x} \qquad (122)$$

$$\overline{x}^T \triangleq [\overline{X}_1^T, \overline{X}_2^T, \cdots, \overline{X}_m^T], \overline{X}_i \in C^{\sigma_i \times 1} \qquad (123)$$

$$\overline{X}_i = \overline{A}_{ii}\overline{X}_i + \overline{u}_i, i=1,2,\cdots,m \qquad (124)$$

The above analysis shows that the regular static state feedback (112) can transform the Luenberger second controllable canonical form into a new controllable system $\sum(A_{c3}, B_{c3}, C_{c2})$ with only $m$ independent chains of integrators with each order of $\sigma_i$, and this system is a typical tree-like linear system, which can be used to analyze the subsequent decoupling process.

**Definition 2.32.** For a controllable linear time-invariant system $\sum(A, B, C)$, a controllable system $\sum(A_{c3}, B_{c3}, C_{c2})$ is called a third controllable standard system (to distinguish it from the Luenberger first and second controllable canonical forms) specially used for decoupling analysis in the latter part.

Obviously, the controllable system $\sum(A_{c3}, B_{c3}, C_{c2})$ is composed of $m$ mutually independent controllable subsystems (i.e., $m$ mutually independent paths) and each subsystem as shown in equation (124) is independently controlled by only one control input, so that there is one and only one independent static state feedback for each subsystem (that is a SISO controllable linear subsystem) to assign the closed-loop poles of this subsystem.

In the following part, it is assumed that the given controllable linear system in Morgan's problem is a tree-like controllable linear system. If it is not a tree-like system, the system is first transformed into a Luenberger second controllable canonical form by a coordinate transformation, then the Luenberger second controllable canonical form is transformed into a third controllable standard system by a static state feedback, and so a tree-like controllable linear system for decoupling is obtained.

If a third controllable standard system $\sum(A_{c3}, B_{c3}, C_{c2})$ can be decoupled into a regular decoupled third controllable standard system $\sum(\hat{A}_{c3}, \hat{B}_{c3}, C_{c2})$ though a static state feedback, it is assumed that the set of the remaining first-order branches except the controllable and observable first-order branches of the decoupled system $\sum(\hat{A}_{c3}, \hat{B}_{c3}, C_{c2})$ is $\overline{x}'$, then an equivalent decoupled third controllable standard system which is used to analyze the pole assignment of the decoupled system can be derived as follows.

$$\hat{x}^T \triangleq \begin{bmatrix} Y^T, (\overline{x}')^T \end{bmatrix} \qquad (125)$$

$$Y^T \triangleq \begin{bmatrix} Y_1^T, Y_2^T, \cdots, Y_p^T \end{bmatrix} \qquad (126)$$

$$Y_i^T \triangleq \begin{bmatrix} y_i, y_i^{(1)}, \cdots, y_i^{(d_i-1)} \end{bmatrix}, i=1,2,\cdots,p \qquad (127)$$

$$Y_i = P_i\overline{x}, i=1,2,\cdots,p \qquad (128)$$

$$P_i \triangleq \begin{bmatrix} C_{c2;i} \\ C_{c2;i}\hat{A}_{c3} \\ \vdots \\ C_{c2;i}\hat{A}_{c3}^{d_i-1} \end{bmatrix}, i=1,2,\cdots,p \qquad (129)$$

$$P \triangleq \begin{bmatrix} P_1^T, P_2^T, \cdots, P_p^T \end{bmatrix}^T \triangleq \begin{bmatrix} \overline{P}_1 & \overline{P}_2 \end{bmatrix} \qquad (130)$$

$$n_{co} \triangleq \sum_{i=1}^p d_i \qquad (131)$$

$$rank P = n_{co} \qquad (132)$$

$$\hat{x} = T\overline{x} \triangleq \begin{bmatrix} \overline{P}_1 & \overline{P}_2 \\ 0 & I \end{bmatrix} \overline{x} \qquad (133)$$

$$\det T \neq 0 \qquad (134)$$

$$\begin{cases} \dot{\hat{x}} = T\hat{A}_{c3}T^{-1}\hat{x} + T\hat{B}_{c3}\upsilon \triangleq \tilde{A}_{c3}\hat{x} + \tilde{B}_{c3}\upsilon \\ y = C_{c2}T^{-1}\hat{x} \triangleq \tilde{C}_{c2}\hat{x} \end{cases} \qquad (135)$$

$$y_i^{(d_i)} = \upsilon_i, i=1,2,\cdots,p \qquad (136)$$

$$\dot{Y}_i = \tilde{A}_{ii}Y_i + b_i\upsilon_i, i=1,2,\cdots,p \qquad (137)$$

$$\tilde{A}_{ii} = \begin{bmatrix} 0 & 1 & 0 & \cdots & 0 \\ 0 & \ddots & \ddots & \ddots & \vdots \\ \vdots & \ddots & \ddots & \ddots & 0 \\ 0 & \cdots & 0 & 0 & 1 \\ 0 & 0 & 0 & 0 & 0 \end{bmatrix} \qquad (138)$$

$$b_i = \begin{bmatrix} 0 & 0 & \cdots & 0 & 1 \end{bmatrix}^T \qquad (139)$$

where $C_{c2;i}$ is the *ith* row of $p \times n$ matrix $C_{c2}$.

The equations (137) represent that the controllable and observable subsystem of the equivalent decoupled third controllable standard system $\sum(\tilde{A}_{c3}, \tilde{B}_{c3}, \tilde{C}_{c2})$ is only composed of $p$ independent chains of integrators with each order of $d_i$ (i.e., $p$ mutually independent forward paths). Moreover, the controllable and unobservable subsystem of $\sum(\tilde{A}_{c3}, \tilde{B}_{c3}, \tilde{C}_{c2})$ is composed of at most $p$ mutually independent paths, but it cannot be guaranteed that each first-order branch in the controllable and unobservable subsystem of $\sum(\tilde{A}_{c3}, \tilde{B}_{c3}, \tilde{C}_{c2})$ is connected to one and only one control input by path. So it cannot be guaranteed that the controllable subsystem of $\sum(\tilde{A}_{c3}, \tilde{B}_{c3}, \tilde{C}_{c2})$ must be composed of $p$ mutually independent controllable subsystems.

In the following part, it is assumed that the controllable subsystem of a regular decoupled linear system used to analyze this subsystem's pole assignment is an equivalent decoupled tree-like controllable subsystem (as shown in (29)). If it is not, the subsystem is first transformed into a decoupled third controllable standard subsystem, then the decoupled third controllable standard subsystem is transformed into an equivalent decoupled third controllable standard subsystem as above (as shown in (135)), and so a regular equivalent decoupled tree-like controllable subsystem used to analyze its pole assignment is obtained.

**Theorem 2.33.** A regular decoupled third controllable standard system $\Sigma(\hat{A}_{c3}, \hat{B}_{c3}, C_{c2})$ can be transformed into a regular equivalent decoupled third controllable standard system $\Sigma(\tilde{A}_{c3}, \tilde{B}_{c3}, \tilde{C}_{c2})$ by a coordinate transformation as shown in (133).

**Proposition 2.34.** For a MIMO controllable linear time-invariant system $\Sigma(A, B, C)$, there exists a regular static state feedback as shown in (141) to assign all closed-loop poles of the system as follows (Wolovich, 1974; Wonham, 1974; Pang, 1992).

$$\begin{cases} \dot{x} = Ax + Bu \\ y = Cx \end{cases} \quad (140)$$

$$u = Kx + \bar{u} \quad (141)$$

$$\dot{x} = (A + BK)x + B\bar{u} \triangleq \bar{A}x + B\bar{u} \quad (142)$$

It is assumed that the designated closed-loop poles of the system are $p_1, p_2, \cdots, p_n$ and the characteristic polynomial of $n \times n$ matrix $\bar{A}$ is $\bar{\Delta}(s)$, the following conclusions can be get by further derivation.

$$\bar{\Delta}(s) = \prod_{j=1}^{n}(s - p_j) \triangleq \sum_{j=0}^{n} \alpha_j s^j, \alpha_n = 1 \quad (143)$$

$$\bar{\Delta}(\bar{A}) = \sum_{j=0}^{n} \alpha_j \bar{A}^j = \prod_{j=1}^{n}(\bar{A} - p_j I) = 0 \quad (144)$$

$$y_i^{(j)} = C_i \bar{A}^j x, 0 \leq j < d_i \quad (145)$$

$$y_i^{(d_i)} = C_i \bar{A}^{d_i} x + C_i \bar{A}^{d_i - 1} B \bar{u} \quad (146)$$

$$y_i^{(j)} = C_i \bar{A}^j x + \sum_{k=d_i}^{j} C_i \bar{A}^{k-1} B \bar{u}^{(k-d_i)}, d_i \leq j \leq n \quad (147)$$

$$\sum_{j=0}^{n} \alpha_j y_i^{(j)} = \sum_{j=0}^{n} \alpha_j C_i \bar{A}^j x + \sum_{j=d_i}^{n} (\sum_{k=d_i}^{j} \alpha_k C_i \bar{A}^{k-1} B \bar{u}^{(k-d_i)}) \quad (148)$$

$$\sum_{j=0}^{n} \alpha_j s^j y_i(s) = C_i (\sum_{j=0}^{n} \alpha_j \bar{A}^j) x(s) + \sum_{j=d_i}^{n} \sum_{k=d_i}^{j} \alpha_k C_i \bar{A}^{k-1} B s^{k-d_i} \bar{u}(s) \quad (149)$$

$$\bar{\Delta}(s) y_i(s) = C_i \bar{\Delta}(\bar{A}) x(s) + \sum_{j=d_i}^{n} \sum_{k=d_i}^{j} \alpha_k C_i \bar{A}^{k-1} B s^{k-d_i} \bar{u}(s) \quad (150)$$

$$m_i(s) \triangleq \sum_{j=d_i}^{n} \sum_{k=d_i}^{j} \alpha_k C_i \bar{A}^{k-1} B s^{k-d_i} \quad (151)$$

$$m_i(s) \triangleq \begin{bmatrix} m_{i1}(s) & m_{i2}(s) & \cdots & m_{ip}(s) \end{bmatrix} \quad (152)$$

$$m(s) \triangleq \begin{bmatrix} m_1(s)^T & m_2(s)^T & \cdots & m_p(s)^T \end{bmatrix}^T \quad (153)$$

$$\bar{\Delta}(s) y_i(s) = C_i \bar{\Delta}(\bar{A}) x(s) + m_i(s) \bar{u}(s) = m_i(s) \bar{u}(s) \quad (154)$$

$$\bar{\Delta}(s) y(s) = m(s) \bar{u}(s) \quad (155)$$

$$T_F \triangleq y(s) / \bar{u}(s) = m(s) / \bar{\Delta}(s) \quad (156)$$

**Proposition 2.35.** For a MIMO linear time-invariant system $\Sigma(A, B, C)$, a static state feedback can only be used to assign the closed-loop poles of the controllable subsystem but not the closed-loop poles of the uncontrollable subsystem.

## 3. The first necessary condition for decoupling system

According to Proposition 2.12, the transfer function $T(s)$ is right invertible. It will first be proved below that right reversibility of the transfer function matrix $T(s)$ means that there exists at least one general $p$-path in the corresponding signal flow graph.

**Theorem 3.1.** For a regular controllable system $\Sigma(A, B, C)$, if its transfer function $T(s)$ is invertible, then there exists at least one general $p$-path in its signal flow graph.

**Proof.** It is first assumed that there are mutually node-joint forward paths from $p$ control inputs to two different outputs such as $y_i$ and $y_j$ though a node (assumed to be $V_0$ which may be a state variable or control input). A row vector composed of the transfer functions of all the forward paths from the $p$ control inputs to the output $y_i$ through the node $V_0$ is $\Delta T_i$, and another row vector composed of the transfer functions of all the forward paths from $p$ control inputs to the output $y_j$ through the same node $V_0$ is $\Delta T_j$.

$$T(s) \triangleq \begin{bmatrix} T_1^T, T_2^T, \cdots T_p^T \end{bmatrix}^T \quad (157)$$

It is obvious that $\Delta T_i$ and $\Delta T_j$ are not equal to zero and have a linear relationship with each other. Row transformations are performed on $T(s)$ to eliminate $\Delta T_i$ or $\Delta T_j$ so that there are no transfer functions of mutually node-joint forward paths from $p$ control inputs to two outputs $y_i$ and $y_j$ though $V_0$ and the invertibility of $T(s)$ is maintained.

$$\Delta T_j = \lambda \Delta T_i, \lambda \neq 0, \lambda \in R(s) \quad (158)$$

$$\Delta T_i = \sum_{k=1}^{p} \gamma_k T_k \quad (159)$$

$$\begin{cases} \hat{T}_i \triangleq T_i - \Delta T_i \\ \hat{T}_j \triangleq T_j - \Delta T_j = T_j - \lambda \Delta T_i \end{cases} \quad (160)$$

Since the $p$ row vectors in $T(s)$ are mutually linearly independent, so

$$\begin{cases} T_j \neq \lambda T_i \\ \hat{T}_j \neq \lambda \hat{T}_i \end{cases} \quad (161)$$

The above equations mean that $\hat{T}_i$ and $\hat{T}_j$ cannot be zero at the same time.

1) One of $\hat{T}_i$ and $\hat{T}_j$ is zero, it is no problem to assume that $\hat{T}_i$ is zero (that is, $\hat{T}_i = 0, \hat{T}_j \neq 0$). Perform the following row transformations as follows.

$$\begin{cases} T_i = \Delta T_i \\ \hat{T}_j = T_j - \Delta T_j = T_j - \lambda \Delta T_i = T_j - \lambda T_i \end{cases} \quad (162)$$

Equations (162) mean that in the matrix $T(s)$, all rows except row $j$ are kept unchanged and a row transformation (denoted as $I(j, i; -\lambda)$, that is performing a row transformation of row $i$ multiplied by a factor of $(-\lambda)$ to superimpose on row $j$) is only performed on row $j$. The new matrix $T(s)$ after the row transformation is also invertible.

2) $\hat{T}_i \neq 0, \hat{T}_j \neq 0$

$$\begin{cases} \hat{T}_i = (1-\gamma_i)T_i - \sum_{k \neq i}^{p} \gamma_k T_k \\ \hat{T}_j = (1-\lambda\gamma_j)T_j - \lambda \sum_{k \neq j}^{p} \gamma_k T_k \end{cases} \quad (163)$$

(a) $1-\gamma_i \neq 0$ or $1-\lambda\gamma_j \neq 0$, it is no problem to assume that $1-\gamma_i \neq 0$.

Equations (163) mean that in the matrix $T(s)$, all rows except row $i$ are kept unchanged and a series of row transformations in turn are performed on row $i$ as follows.

$$\begin{cases} I(i; 1-\gamma_i) \\ I(i,i'; -\gamma_{i'}), i' \neq i, i' = 1, 2, \cdots, p \end{cases} \quad (164)$$

where $I(i; 1-\gamma_i)$ is performing a row transformation of row $i$ multiplied by a factor of $(1-\gamma_i)$.

Similarly, the new matrix $T(s)$ after above row transformations is also invertible.

(b) $1-\gamma_i = 0$ and $1-\lambda\gamma_j = 0$, so $\gamma_i \neq 0, \gamma_j \neq 0$.

$$\hat{T}_i = -\gamma_j T_j - \sum_{\substack{k \neq i \\ k \neq j}}^{p} \gamma_k T_k = -\frac{1}{\lambda} T_j - \sum_{\substack{k \neq i \\ k \neq j}}^{p} \gamma_k T_k \quad (165)$$

Equation (165) means that in the matrix $T(s)$, all rows except row $j$ are kept unchanged and a series of row transformations are performed on row $j$ as follows.

$$\begin{cases} I(j; -1/\lambda) \\ I(j,i'; -\gamma_{i'}), i' \neq i, i' \neq j, i' = 1, 2, \cdots, p \end{cases} \quad (166)$$

The row $j$ of the new matrix $T(s)$ after above row transformations is transformed into the row vector $\hat{T}_i$.

In turn, a series of row transformations are performed on the matrix $T(s)$ to eliminate the influence of the transfer functions of the mutually node-joint forward paths from $p$ control inputs to different outputs though the same node and the invertibility of $T(s)$ is maintained until there are no transfer functions of the mutually node-joint forward paths from $p$ control inputs to different outputs though the same node in $T(s)$.

It is then assumed that there are forward paths from the two different control inputs such as $u_i$ and $u_j$ to one output such as $y_v$. Then, a series of the column transformations (that are similar to the above row transformations) are performed on the matrix $T(s)$ to eliminate the transfer functions of the forward paths from $u_i$ or $u_j$ to output $y_v$ and the invertibility of $T(s)$ is maintained. In turn, a series of column transformations are performed on the matrix $T(s)$ to eliminate the transfer functions of the forward paths from the different control inputs to one output and the invertibility of $T(s)$ is maintained until there are no transfer functions of the forward paths from the different control inputs to the same output in $T(s)$.

Therefore, the matrix $T(s)$ is transformed into a new matrix whose each row (or column) has one and only one non-zero item. It means that there exists at least one general $p$-path in the corresponding signal flow graph, and so Theorem 3.1 holds.

□

The following theorems are easy derived from Theorem 3.1.

**Theorem 3.2.** For a nonregular controllable system $\sum(A, B, C)$, if its transfer function $T(s)$ is right invertible, then there exists at least one general $p$-path in its signal flow graph.

**Theorem 3.3.** For a nonregular controllable system $\sum(A, B, C)$, if $\sum(A, B, C)$ can be decoupled, then there exists at least one general $p$-path in the signal flow graph corresponding to its decoupled system.

**Definition 3.4.** The general relative order (denoted as $r_i$) for the output (such as $y_i$) of a general path is defined as the order or length of this general path. The master (control) inputs are defined as the control inputs of all $p$ general paths in one general $p$-path.

The decoupling process of the system can be equivalent to the process of condensing the forward paths with the same output as one general path to this general path so that the control inputs that need to be compensated are obtained. Once the system has been decoupled, the controllable and observable subsystem of the decoupled system is equivalent to condensing into an equivalent decoupled system.

Moreover, each general path is required to be the shortest path corresponding to its control input and output; otherwise, it is assumed that a general path such as $path(u_j, y_i)$ is not the shortest path such as $path(u_j, y_i)_{min}$, then the derivative of the output $y_i$ (whose derivative order is its general relative order) along this general path is obtained, and there must be a derivative term of the control input $u_j$ in the obtained control input item for the derivative of the output $y_i$ so that the system can never be decoupled, no matter how many integrator-strings are used to compensate before the input $u_j$ of the general path. Therefore, the decoupling framework is defined as follows.

**Definition 3.5.** For a nonregular controllable system $\sum(A, B, C)$, if its transfer function $T(s)$ is right invertible, a decoupling framework (if exists) is defined as one general $p$-path in the signal flow graph, each of which is the shortest path corresponding to its control input and output.

**Theorem 3.6.** For a nonregular controllable system $\sum(A, B, C)$, the first necessary condition for decoupling this system is that there exists at least one decoupling framework in the corresponding signal flow graph.

It should be noted here that there may be more than one decoupling frameworks, so all decoupling frameworks must be exhausted in the decoupling process.

**Theorem 3.7.** For a nonregular controllable system $\sum(A, B, C)$, if $\sum(A, B, C)$ can be decoupled, then there exists at least one decoupling framework in the signal flow graph corresponding to its decoupled system.

**Theorem 3.8.** For a nonregular controllable system $\sum(A, B, C)$, if $\sum(A, B, C)$ can be decoupled, then there exists one and only one decoupling framework in the signal flow graph corresponding to its equivalent decoupled tree-like system.

## 4. The second necessary condition for decoupling system

**Definition 4.1.** For a nonregular controllable system $\sum(A, B, C)$, if there exists one decoupling framework in its signal flow graph, a set of $p$ general relative orders is $\{r_i\}_p$ for this decoupling framework and a set of $p$ relative orders is $\{d_i\}_p$, then the following inequality holds.
$$r_i \geq d_i, \forall i \quad (167)$$

$$E_{de} = \begin{bmatrix} \sum_{j=d_1}^{r_1-1} C_1 A^{j-1} B s^{r_1-j} \\ \sum_{j=d_2}^{r_2-1} C_2 A^{j-1} B s^{r_2-j} \\ \vdots \\ \sum_{j=d_p}^{r_p-1} C_p A^{j-1} B s^{r_p-j} \end{bmatrix} \quad (168)$$

The matrix $E_{de}$ is called an independent decoupling matrix corresponding to each decoupling framework. Through the independent decoupling matrix, it can be initially obtained that each non-zero column of the matrix $E_{de}$ corresponding to the control input needs to be compensated and the highest order of $s$ among all elements in this column is its compensation order. To determine the control inputs that need to be compensated and their required compensation orders, a complicated process is also required.

It is assumed that there are $p_1$ integrator-strings which are node-disjoint with the general $p$-path in a decoupling framework can be selected to compensate control inputs.

**Definition 4.2.** For a nonregular controllable system $\Sigma(A,B,C)$, if there exists one decoupling framework in the corresponding signal flow graph and there are $p_1$ integrator-strings which are node-disjoint with the general $p$-path, a general decoupling framework is defined as one general $(p+p_1)$-path composed of the $p_1$ integrator-strings selected and the general $p$-path, each of which is the shortest path corresponding to its control input and output.

It has been assumed that all integrator-strings used to compensate are already selected into the general decoupling framework.

The control input of each integrator-string selected is called a general master (control) input and the remaining control inputs except all general master inputs and master inputs in $u$ are called auxiliary (control) inputs. Then the output of each integrator-string selected is called a general output and a column vector (or set) composed of all general outputs is denoted as $y'$.

$$y' \triangleq \begin{bmatrix} C_{p+1}^T & C_{p+2}^T & \cdots & C_{\hat{p}}^T \end{bmatrix}^T x, \hat{p} \triangleq p + p_1 \quad (169)$$

$$\hat{y} \triangleq \begin{bmatrix} y \\ y' \end{bmatrix} = \begin{bmatrix} C_1^T & C_2^T & \cdots & C_{\hat{p}}^T \end{bmatrix}^T x \triangleq \hat{C}x \quad (170)$$

$$\hat{E}_{de} = \begin{bmatrix} \sum_{j=d_1}^{r_1-1} C_1 A^{j-1} B s^{r_1-j} \\ \sum_{j=d_2}^{r_2-1} C_2 A^{j-1} B s^{r_2-j} \\ \vdots \\ \sum_{j=d_{\hat{p}}}^{r_{\hat{p}}-1} C_{\hat{p}} A^{j-1} B s^{r_{\hat{p}}-j} \end{bmatrix} \quad (171)$$

The matrix $\hat{E}_{de}$ is called a general independent decoupling matrix corresponding to each general decoupling framework. On the assumption that all auxiliary inputs can be self-compensated, their corresponding columns in the matrix $\hat{E}_{de}$ can be deleted first. The auxiliary inputs (called first-type auxiliary (control) inputs) which need self-compensation in the process of the compensations for the below master inputs (or general master inputs) are preset to zero. Finally, the master inputs (or general master inputs) that need to be compensated and their required compensation orders are determined by use of the matrix $\hat{E}_{de}$.

It is no problem to assume that all control inputs corresponding to all columns in the matrix $\hat{E}_{de}$ are master inputs or general master inputs.

It should be noted that there may be more than one selections of integrator-strings in each decoupling framework, so all possible selections of integrator-strings must be exhausted in the decoupling process.

The output of an integrator-string may be selected in order as follows: 1) the state variable (the output node of whose corresponding general first-order branch is not connected to other first-order branches) in an integrator-string is firstly selected; 2) the state variable (the output node of whose corresponding general first-order branch is also the output node of at least one other general first-order branch) in an integrator-string is secondly selected; and 3) the other state variable in an integrator-string is selected.

If the order of an integrator-string selected is greater than the required compensation order, there are more than one compensation options: 1) all first-order branches in this integrator-string are used to compensate, then these first-order branches are controllable and observable in the decoupled system once the system can be decoupled, so that they can be assigned their closed-loop poles; and 2) non-all first-order branches in this integrator-string are used to compensate, then the first-order branches which are not used to compensate in this integrator-string may be controllable and unobservable in the decoupled system once the system can be decoupled, so that whether they can be assigned their closed-loop poles depends on whether the decoupled controllable subsystem is composed of $p$ mutually independent controllable subsystems (see Section 8).

**Definition 4.3.** For a nonregular controllable system $\Sigma(A,B,C)$, if there are one decoupling framework and its one general decoupling framework, the matrix $\hat{E}_{de}$ may be decomposed into a series of independent sub-matrices (called general independent decoupling sub-matrices denoted as $\hat{E}_{de}^k$), in each of which, each row (or column) vector is non-zero and the remaining elements (except all elements in each row (or column) vector of this independent sub-matrix) in corresponding row (or column) vector of the matrix $\hat{E}_{de}$ are zero. Obviously, these decomposed general independent decoupling sub-matrices are unique and mutually independent, and can determine the control inputs that need to be compensated and their required compensation orders.

If one general independent decoupling sub-matrix is a column vector, its corresponding control input (called master decoupling input) needs to be compensated by an integrator-string and the highest order of $s$ among all elements in this column vector is its required compensation order. This compensation method that directly compensates the integrator-string for a control input is called the first-type compensation by chain of integrators (that is an independent compensation method).

If one general independent decoupling sub-matrix has more than one columns, there may be many options to determine the control inputs that need to be compensated, their required compensation orders and compensation expressions. So all possible selections must be exhausted in the decoupling process.

1) The control input corresponding to each column in one general independent decoupling sub-matrix is independently compensated as same as above, but this is not an optimal compensation method.

2) It is assumed that the column number in one sub-matrix such as $\hat{E}_{de}^k$ is $\rho_k$. The control inputs corresponding to $\rho_{k1}$ (there are many choices while $\rho_{k1}$ changes from 1 to $\rho_k$) columns are selected as master decoupling inputs (denoted as $\tilde{u}^k$), and the remaining control inputs (i.e., $\rho_{k2}(=\rho_k-\rho_{k1})$ inputs) except master decoupling inputs are auxiliary decoupling inputs (denoted as $\bar{u}^k$). Each master decoupling input which needs to be compensated and its required compensation order can be determined immediately by the corresponding column vector. Since the column vector is already non-zero, its corresponding control input needs to be compensated by an integrator-string and the highest order of $s$ among all elements in this column vector is its required compensation order. It is called the second-type compensation by chain of integrators (that is a block compensation method). The undetermined coefficient method is used to determine the compensation expressions of all master decoupling inputs in the following.

$$\tilde{u}^k \triangleq \left[\tilde{u}_1^k, \tilde{u}_2^k, \cdots, \tilde{u}_{\rho_{k1}}^k\right]^T \quad (172)$$

$$\bar{u}^k \triangleq \left[\bar{u}_1^k, \bar{u}_2^k, \cdots, \bar{u}_{\rho_{k2}}^k\right]^T \quad (173)$$

$$u^k \triangleq \tilde{u}^k \cup \bar{u}^k \quad (174)$$

It is assumed that the highest order of $s$ for one master decoupling input $\tilde{u}_i^k$ in the row $j$ of $\hat{E}_{de}^k$ is $\delta_{ij}^k$ and the corresponding row in the matrix $\hat{E}_{de}$ is row $k_j$.

$$\hat{E}_{de;k_j} \triangleq \sum_{\mu=d_{k_j}}^{r_{k_j}-1} C_{k_j} A^{\mu-1} B s^{r_{k_j}-\mu} \quad (175)$$

$$\bar{E}_{ij}^k \triangleq \sum_{\mu=d_{k_j}}^{r_{k_j}-\delta_{ij}^k} C_{k_j} A^{\mu-1} B s^{r_{k_j}-\delta_{ij}^k-\mu} \quad (176)$$

$$\bar{D}_{ij}^k \triangleq C_{k_j} A^{r_{k_j}-\delta_{ij}^k} \quad (177)$$

On the right side of equation (175), all items of the master decoupling inputs and auxiliary decoupling inputs (belonging to this sub-matrix $\hat{E}_{de}^k$) are retained, but the remaining items are deleted so that a new row vector (denoted as $\hat{E}_{ij}^k$) is get.

On the right side of equation (176), all items of the master decoupling inputs (belonging to this sub-matrix $\hat{E}_{de}^k$) and the control inputs (belonging to other general independent decoupling sub-matrices) are deleted, but all items of the auxiliary decoupling inputs (belonging to this sub-matrix $\hat{E}_{de}^k$) are retained so that a new row vector (denoted as $\tilde{E}_{ij}^k$) is obtained.

On the right side of equation (177), each non-zero item whose corresponding state variable is connected with one control input (belonging to this sub-matrix $\hat{E}_{de}^k$) by path is retained and the remaining items are set to zero so that a new row vector (denoted as $\tilde{D}_{ij}^k$) is get.

It is assumed that the undetermined coefficients are $\alpha_{ij}^k \triangleq diag[\alpha_{ij1}^k, \alpha_{ij2}^k, \cdots, \alpha_{ij\rho_{k2}}^k]$ (whose corresponding column is set to zero if the column in $\tilde{E}_{ij}^k$ is zero) and $\beta_{ij}^k \triangleq diag[\beta_{ij1}^k, \beta_{ij2}^k, \cdots, \beta_{ijn}^k]$ (whose corresponding column is set to zero if the column in $\tilde{D}_{ij}^k$ which is corresponding to one state variable is zero).

Then let $j$ change from 1 to $\eta_k$ (that is defined as the row number of the matrix $\hat{E}_{de}^k$), all matrices $\tilde{D}_{ij}^k$ and $\tilde{E}_{ij}^k$ corresponding to $\tilde{u}_i^k$ are obtained in turn.

$$\tilde{u}_i^k = z_i^k + \sum_{j=1}^{\eta_k} \tilde{E}_{ij}^k (\alpha_{ij}^k \bar{u}^k) + \sum_{j=1}^{\eta_k} \tilde{D}_{ij}^k (\beta_{ij}^k x) \quad (178)$$

$$z^k \triangleq \left[z_1^k, z_2^k, \cdots, z_{\rho_{k1}}^k\right]^T \quad (179)$$

$$\tilde{u}^k \triangleq z^k + \tilde{E}_0^k \bar{u}^k + \tilde{D}_0^k x \quad (180)$$

where $z_i^k$ is the general output of an integrator-string (whose general master input is denoted as $u(z_i^k)$). The general path such as $path(\tilde{u}_i^k, y_{k_j})$ is transformed into a new general path (denoted as $path(u(z_i^k), y_{k_j})$) compensated by an integrator-string of $path(u(z_i^k), z_i^k)$. Then, $path(u(z_i^k), y_{k_j})$ is composed of $path(u(z_i^k), z_i^k)$ and $path(\tilde{u}_i^k, y_{k_j})$, and must be the shortest path corresponding to its control input and output; otherwise, this integrator-string cannot be used to compensate input $\tilde{u}_i^k$.

For a certain $i$, while the corresponding undetermined coefficients are set to zero, the input $\tilde{u}_i^k$ will be independently compensated by an integrator-string which is the first-type compensation.

Let $i$ change from 1 to $\rho_{k1}$ (that is the number of master decoupling inputs in $\hat{E}_{de}^k$), the compensation expressions of all master decoupling inputs in $\hat{E}_{de}^k$ are obtained.

Substitute (178) into $\hat{E}_{de}^k u^k$, the derivative term of each state variable is converted into a linear combination of the control inputs or other state variables by use of the state equations of system $\sum(A, B, C)$ and obtain

$$\hat{E}_{de}^k u^k \triangleq \tilde{E}_1^k z^k + \tilde{E}_2^k \bar{u}^k + \tilde{D}_1^k x \quad (181)$$

$$\tilde{E}_3^k \triangleq \begin{bmatrix} \tilde{E}_0^k \\ \tilde{E}_2^k \end{bmatrix} \quad (182)$$

Each column in the matrix $\tilde{E}_3^k$ is corresponding to the auxiliary decoupling input and all elements in the matrix $\tilde{E}_3^k$ are constant or polynomials of $s$.

The purpose of using the second-type compensation by chain of integrators is to reduce the highest order of $s$ among all elements in each column of $\tilde{E}_3^k$ as much as possible by selecting the undetermined coefficients so that the required compensation order corresponding to the column of $\tilde{E}_3^k$ is also reduced, which is beneficial to decouple the system. And the ideal goal is to transform all elements in each column of $\tilde{E}_3^k$ to constant once all undetermined coefficients are determined, so that the corresponding auxiliary decoupling input does not need compensation, which is very beneficial to decouple the system.

It is assumed that the highest order of $s$ among all elements in the column $i$ of $\tilde{E}_3^k$ (denoted as $\tilde{E}_{3;i}^k$) is $\bar{k}_i$.

$$\bar{k}_i \geq 0, i=1,2,\cdots,\rho_{k2} \quad (183)$$

$$\tilde{E}_{3;i}^k \triangleq \sum_{j=0}^{\bar{k}_i} f_{i,j}^k s^j, i=1,2,\cdots,\rho_{k2} \quad (184)$$

where $f_{i,j}^k$ is a $(\rho_{k1}+\eta_k)\times 1$ column vector whose elements are the polynomials about the undetermined coefficients.

If $\bar{k}_i=0$, the auxiliary decoupling input $\bar{u}_i^k$ does not need compensation. It is no problem to assume that all auxiliary decoupling inputs need to be compensated (i.e., $\bar{k}_i \geq 1, \forall i$) in the below analysis.

The equations about the undetermined coefficients are obtained in order to transform all elements in each column of $\tilde{E}_3^k$ to constant.

$$f_{i,j}^k = 0, j = \bar{k}_i, \bar{k}_i - 1, \cdots, 1; i=1,2,\cdots,\rho_{k2} \quad (185)$$

where, if $f_{i,j}^k \equiv 0$, it means that there is no $j$th order item in the column vector $\tilde{E}_{3;i}^k$; otherwise, it means that there are $j$th order items in $\tilde{E}_{3;i}^k$.

For a fixed column $i$, sometimes no matter how the undetermined coefficients are selected, not all equations in (185) can be established. When $j$ changes from $\bar{k}_i$ to one in turn, by selecting the undetermined coefficients, the minimum value of $j$ (denoted as $\min_i^k$, $\min_i^k \geq 1, \forall i$) at which the equations (formed from (185) when $j$ changes from $\bar{k}_i$ to $\min_i^k$) have at least one solution about the undetermined coefficients. That is to say, if $\min_i^k \geq 2$, the equations (formed from (185) when $j$ changes from $\bar{k}_i$ to $\min_i^k$) have at least one solution but the equations (formed from (185) when $j$ changes from $\bar{k}_i$ to $(\min_i^k - 1)$) have no solution; otherwise, $\min_i^k = 1$ holds. It is shown that the order items of $s$ from $\bar{k}_i$ to $\min_i^k$ in the column $i$ of $\tilde{E}_3^k$ can be independently eliminated by selecting the undetermined coefficients, so that the compensation order corresponding to the column vector $\tilde{E}_{3;i}^k$ is reduced to $(\min_i^k - 1)$.

Of course, if there exists a certain $i$, when the equation formed from (185) when $j$ is equal to $\bar{k}_i$ has no solution, it is means that the $\bar{k}_i$th order items of $s$ in the column vector $\tilde{E}_{3;i}^k$ cannot be eliminated and $\min_i^k = \bar{k}_i + 1$ holds. It is obvious that the input $\bar{u}_i^k$ needs to be compensated by a $\bar{k}_i$th order integrator-string.

While $i$ changes from 1 to $\rho_{k2}$, $\min_1^k, \min_2^k, \cdots, \min_{\rho_{k2}}^k$ are obtained.

$$\min^k \triangleq \{\min_i^k\}_{\rho_{k2}} = [\min_1^k, \min_2^k, \cdots, \min_{\rho_{k2}}^k]^T \quad (186)$$

**Definition 4.4.** For a general independent decoupling matrix (or sub-matrix) and one kind of selection of master decoupling inputs, a set of coefficients such as $\tau_1^k$, $\tau_2^k, \cdots, \tau_{\rho_{k2}}^k$ ($\tau_i^k \geq \min_i^k, \forall i$) is defined as a smallest parameter set of the largest equations formed from (185) with at least one solution. That is to say, when $j$ changes from $\bar{k}_i$ to $\tau_i^k$ for each value of $i$, the equations formed from (185) have at least one solution; but, when $j$ changes from $\bar{k}_i$ to $\tau_i^k$ for each value of $i$ except $i'$ ($i' = 1, 2, \cdots, \rho_{k2}$) for which $j$ changes from $\bar{k}_{i'}$ to $(\tau_{i'}^k - 1)$ (if $\tau_{i'}^k \geq 2$), the equations formed from (185) have no solution. And the number of all possible choices for the smallest parameter set of the largest equations formed from (185) with at least one solution is set to $n_{de1}^k$ (called the number of the first-level smallest parameter set for the sub-matrix $\tilde{E}_3^k$).

For each selection of smallest parameter set, the $\tilde{E}_3^k$ will be transformed into a new matrix which the highest order of $s$ among all elements in the column $i$ is reduced to $(\tau_i^k - 1)$. In this new matrix, after the columns or rows (in which all elements are constant) are deleted, if the obtained matrix is non-constant, it is one second level of general independent decoupling sub-matrix (denoted as $\hat{E}_{de2}^k$, and $\hat{E}_{de}^k$ is also called the first level of general independent decoupling sub-matrix) that will be processed below; otherwise, the processing of the general independent decoupling sub-matrix is completed.

Then, the next level of general independent decoupling sub-matrices will be analyzed and processed as same as above first level of general independent decoupling sub-matrices until all the elements in all obtained sub-matrices are constant.

Therefore, the master inputs or general master inputs (which may be the different levels of master decoupling inputs) that need to be compensated by integrator-strings and their required compensation orders are obtained.

Moreover, since the compensation for master inputs (or general master inputs) may add other master inputs (or general master inputs) that need to be compensated, the above results obtained need to be verified and further compensation by integrator-strings may be performed on the new master inputs (or general master inputs) that need to be compensated. This requires a process of judgment and cyclic processing which can be guaranteed to stop after a finite number of cyclic processing. The specific cyclic processing flow is as follows.

1) All auxiliary inputs of the system are preset to zero.

2) The master inputs (or general master inputs) that need to be compensated and their required compensation orders are derived as above from the general independent decoupling matrix.

3) If there are enough integrator-strings with appropriate orders in the general decoupling framework, each master input (or general master input) that needs to be compensated by a chain of integrators with an appropriate order (which may need to be compensated by one or more than one other integrator-strings to increase its order, called re-compensation of this integrator-string) is pre-compensated, and then an updated system, its updated general decoupling framework and the set of its updated general relative orders are obtained, then jump to step 4); otherwise, the original system cannot be decoupled (for the selected general decoupling framework and selected different levels of master decoupling inputs) so that the cyclic processing flow is stopped.

4) The controllable subsystem of the updated system can be transformed into a third controllable standard subsystem, and its controllability and row by row decouplability cannot be changed.

5) The set of updated relative orders is obtained in the third controllable standard subsystem.

6) Though comparing the set of updated general relative orders with the set of updated relative orders, if the two are

equal, it means that the above results do not introduce any new master input (or general master input) that needs to be compensated; otherwise, the above results have introduced at least one new input that needs to be compensated.

7) If at least one new input that needs to be compensated is introduced in the process of compensating the inputs that need to be compensated, the new inputs that need to be compensated and their required compensation orders are derived from the corresponding updated general independent decoupling matrix. Then turn to step 3).

8) If no new input that needs to be compensated is introduced in the process of compensating the inputs that need to be compensated, then the second necessary condition for decoupling system as shown in Theorem 4.5 is met, so that the cyclic processing flow is stopped. Therefore, the compensation expressions (corresponding to the original system) of all master inputs (or general master inputs) that need to be compensated will be finally obtained.

There are a finite number of integrator-strings in the general decoupling framework selected, so the above cyclic processing flow will be stopped after a finite number of cyclic processing, so that either the original system cannot be decoupled (for the selected general decoupling framework and selected different levels of master decoupling inputs), or the second necessary condition for decoupling system is met.

**Theorem 4.5.** For a general decoupling framework and its general independent decoupling matrix, if the corresponding nonregular controllable system $\Sigma(A,B,C)$ can be decoupled, then in at least one general decoupling framework, there are enough integrator-strings with appropriate orders to meet all the master inputs (or general master inputs) that need to be compensated and their required compensation orders derived from the corresponding general independent decoupling matrix (which is the second necessary condition for decoupling a nonregular system).

## 5. A sufficient condition for decoupling system

A pre-decoupling system is obtained after each master input (or general master input) that needs to be compensated is already compensated under the premise that the first-type auxiliary inputs are preset to zero. In each pre-decoupling system, a judgable decoupling system is obtained after the first-type auxiliary inputs have been set to zero, and then a sufficient condition that the judgable decoupling system can be decoupled is given in the following.

It is assumed that a set of master inputs (or general master inputs) that have completed first-type compensation in a general decoupling framework is $U_{01}$, another set of master inputs (or general master inputs) that have completed second-type compensation is $U_{02}$, a set of first-type auxiliary inputs is $\bar{U}_2$, and a set of remaining control inputs except the control inputs of $U_{01}, U_{02}$ and $\bar{U}_2$ in $u$ is $\bar{U}_1$.

$$\begin{cases} \dot{x} = Ax + Bu \\ y = Cx \end{cases} \quad (187)$$

$$\begin{cases} U_{01} = K_{01}x \\ U_{02} = \bar{K}_{02}x + \bar{G}_{02}\bar{U}_1 \\ \bar{U}_2 = 0 \end{cases} \quad (188)$$

$$u^T \triangleq \begin{bmatrix} U_{01}^T & U_{02}^T & \bar{U}_2^T & \bar{U}_1^T \end{bmatrix} \quad (189)$$

$$\bar{U}_1^T \triangleq \begin{bmatrix} u_{v_1} & u_{v_2} & \cdots & u_{v_{p'}} \end{bmatrix}, p' \geq p \quad (190)$$

$$B \triangleq \begin{bmatrix} B_{01} & B_{02} & \bar{B}_2 & \bar{B}_1 \end{bmatrix} \quad (191)$$

Substitute (188) into (187) and obtain as follows.

$$\begin{cases} \dot{x} = A_{D1}x + B_{D1}\bar{U}_1 \\ y = Cx \end{cases} \quad (192)$$

$$\begin{cases} A_{D1} = A + B_{01}K_{01} + B_{02}\bar{K}_{02} \\ B_{D1} = \bar{B}_1 + B_{02}\bar{G}_{02} \end{cases} \quad (193)$$

It is obvious that the system $\Sigma(A_{D1}, B_{D1}, C)$ is one judgable decoupling system corresponding to the pre-decoupling system.

$$B_{D1}^* = \begin{bmatrix} C_1 A_{D1}^{d_1-1} B_{D1} \\ C_2 A_{D1}^{d_2-1} B_{D1} \\ \vdots \\ C_p A_{D1}^{d_p-1} B_{D1} \end{bmatrix} \triangleq \begin{bmatrix} B_1^* & B_2^* & \cdots & B_{p'}^* \end{bmatrix} \quad (194)$$

where $d_1, d_2, \cdots, d_p$ are the relative orders of $\Sigma(A_{D1}, B_{D1}, C)$.

It should be noted here that there is not necessarily a one-to-one correspondence between the general relative orders and relative orders in the system $\Sigma(A_{D1}, B_{D1}, C)$, but there exists a sequence $j_1, j_2, \cdots, j_p$ which is a kind of mathematic permutation of $1, 2, \cdots, p$ so that $d_i = r_{j_i}$, $i = 1, 2, \cdots, p$ (see Example 9.2 in Section 9 for details). Furthermore, there is indeed a one-to-one correspondence between the general relative orders and relative orders in the equivalent decoupled (tree-like) system.

There is a sufficient condition that the system $\Sigma(A, B, C)$ can be decoupled according to Theorem 2.9.

**Theorem 5.1.** For a nonregular Morgan's problem, let

$$\hat{\upsilon} = B_{D1}^* u \quad (195)$$

$$\hat{\upsilon} \triangleq \{\hat{\upsilon}_i\}_p \triangleq [\hat{\upsilon}_1, \hat{\upsilon}_2, \cdots, \hat{\upsilon}_p]^T \quad (196)$$

Then there exists a static state feedback which decouple the system $\Sigma(A_{D1}, B_{D1}, C)$ if $\hat{\upsilon}_1, \hat{\upsilon}_2, \cdots, \hat{\upsilon}_p$ are linearly independent that is equivalent to $rank B_{D1}^* = p$ (which is a sufficient condition for decoupling a nonregular system).

If $rank B_{D1}^* = p$, there is no problem to assume that the first $p$ columns of the matrix $B_{D1}^*$ are linearly independent, and so a set of master inputs of $\Sigma(A_{D1}, B_{D1}, C)$ (i.e., the control inputs of its decoupled system) which is composed of the corresponding first $p$ rows of $\bar{U}_1$ is set to $U_1$.

$$\bar{B}_{D1}^* \triangleq [B_1^* \quad B_2^* \quad \cdots \quad B_p^*] \quad (197)$$

$$\det \bar{B}_{D1}^* \neq 0 \quad (198)$$

$$\bar{U}_1 \triangleq \begin{bmatrix} U_1 \\ \tilde{U}_1 \end{bmatrix} \quad (199)$$

$$U_2 \triangleq \begin{bmatrix} \bar{U}_2 \\ \tilde{U}_1 \end{bmatrix} \quad (200)$$

$$\bar{G}_{02} \triangleq \begin{bmatrix} G_{02} & \tilde{G}_{02} \end{bmatrix} \tag{201}$$

where $\tilde{U}_1$ is a set of remaining auxiliary inputs (called second-type auxiliary (control) inputs) except the first-type auxiliary inputs, and $U_2$ is all auxiliary inputs in the pre-decoupling system which need to be self-compensated.

$$\begin{cases} U_{01} = K_{01}x \\ U_{02} = \bar{K}_{02}x + \bar{G}_{02}\bar{U}_1 = \bar{K}_{02}x + G_{02}U_1 + \tilde{G}_{02}\tilde{U}_1 \end{cases} \tag{202}$$

$$B \triangleq \begin{bmatrix} B_{01} & B_{02} & B_2 & B_1 \end{bmatrix} \tag{203}$$

$$u^T \triangleq \begin{bmatrix} U_{01}^T & U_{02}^T & U_2^T & U_1^T \end{bmatrix} \tag{204}$$

Substitute (202) into (187) and obtain as follows.

$$\begin{cases} \dot{x} = A_{D2}x + B_{D21}U_1 + B_{D22}U_2 \\ y = Cx \end{cases} \tag{205}$$

$$\begin{cases} A_{D2} = A + B_{01}K_{01} + B_{02}\bar{K}_{02} \\ B_{D21} = B_1 + B_{02}G_{02} \\ B_{D22} = B_2 + \begin{bmatrix} 0 & B_{02}\tilde{G}_{02} \end{bmatrix} \end{cases} \tag{206}$$

According to Definition 2.24, the self-compensation of all auxiliary inputs $U_2$ is shown in the following.

$$U_2 = \begin{bmatrix} \bar{U}_2 \\ \tilde{U}_1 \end{bmatrix} = F_{22}x \triangleq \begin{bmatrix} K_{21} \\ K_{22} \end{bmatrix} x \tag{207}$$

$$\tilde{U}_1 = K_{22}x \tag{208}$$

$$\begin{cases} U_{01} = K_{01}x \\ U_{02} = (\bar{K}_{02} + \tilde{G}_{02}K_{22})x + G_{02}U_1 \triangleq K_{02}x + G_{02}U_1 \\ U_2 = F_{22}x \end{cases} \tag{209}$$

Substitute (209) into (205) and obtain

$$\begin{cases} \dot{x} = A_{D3}x + B_{D3}U_1 \\ y = Cx \end{cases} \tag{210}$$

$$\begin{cases} A_{D3} = A_{D2} + B_{D22}F_{22} \\ B_{D3} = B_{D21} \end{cases} \tag{211}$$

Therefore, the system $\Sigma(A_{D3}, B_{D3}, C)$ is a regular decoupleable system that ensures the stability of the uncontrollable subsystem of its decoupled system by assigning the closed-loop poles of this uncontrollable subsystem during the self-compensation process of all auxiliary inputs, as long as a given linear system is controllable.

$$B_{D3}^* = \begin{bmatrix} C_1 A_{D3}^{d_1-1} B_{D3} \\ C_2 A_{D3}^{d_2-1} B_{D3} \\ \vdots \\ C_p A_{D3}^{d_p-1} B_{D3} \end{bmatrix} = \bar{B}_{D1}^* \tag{212}$$

$$A_{D3}^* = \begin{bmatrix} C_1 A_{D3}^{d_1} \\ C_2 A_{D3}^{d_2} \\ \vdots \\ C_p A_{D3}^{d_p} \end{bmatrix} \tag{213}$$

$$\upsilon = A_{D3}^* x + B_{D3}^* U_1 = A_{D3}^* x + \bar{B}_{D1}^* U_1 \tag{214}$$

$$U_1 = F_1 x + G_1 \upsilon \tag{215}$$

$$\begin{cases} F_1 = -(\bar{B}_{D1}^*)^{-1} A_{D3}^* \\ G_1 = (\bar{B}_{D1}^*)^{-1} \end{cases} \tag{216}$$

**Theorem 5.2.** For a regular system $\Sigma(A_{D3}, B_{D3}, C)$ as shown above, there exists a static state feedback which decouple this system if and only if $\det B_{D3}^* = \det \bar{B}_{D1}^* \neq 0$.

## 6. The sufficient and necessary condition for decouping system

In the above part, the first and second necessary conditions, and a sufficient condition that the system can be decoupled are obtained in sequence, and these three conditions are nested sequentially. Moreover, as long as all possible situations that may meet the first and second necessary conditions for decoupling system are exhausted, it can be determined by the sufficient condition (as shown in Theorem 5.1) whether the system can be decoupled.

**Theorem 6.1.** The sufficient and necessary condition for decoupling system is that the system can be decoupled as long as one situation meets these three conditions for decoupling system at the same time, i.e., there exists at least one system $\Sigma(A_{D1}, B_{D1}, C)$ which can be decoupled; otherwise, the system cannot be decoupled.

## 7. Static state feedback for decoupling system

Finally, according to the above decoupling process, it is easy to obtain the static state feedback for decoupling system. Combine (209) and (215) to obtain the nonregular static state feedback $(F, G)$ which can decouple the system $\Sigma(A, B, C)$ as follows.

$$u = \begin{bmatrix} U_{01} \\ U_{02} \\ U_2 \\ U_1 \end{bmatrix} = Fx + G\upsilon \tag{217}$$

$$F = \begin{bmatrix} K_{01} \\ K_{02} + G_{02}F_1 \\ F_{22} \\ F_1 \end{bmatrix}, G = \begin{bmatrix} 0 \\ G_{02}G_1 \\ 0 \\ G_1 \end{bmatrix} \tag{218}$$

$$\det G_1 \neq 0 \tag{219}$$

$$rank G = p \tag{220}$$

The closed-loop transfer function is

$$T_{FG}(s) = diag[s^{-d_1}, s^{-d_2}, \cdots, s^{-d_p}] \tag{221}$$

## 8. The pole assignment for a decoupled system

In this section, the pole assignment for the controllable subsystem of a regular equivalent decoupled tree-like linear system is only analyzed, and so there is no problem to assume that the regular equivalent decoupled tree-like linear system is controllable.

For a regular equivalent decoupled tree-like controllable linear system as shown in (29), in order to assign its closed-loop poles while ensuring its decoupled state, it can be seen from the equation (154): after the closed-loop poles are assigned, it must be guaranteed that there is one and only one control input connected with each output (such as $y_i$) by path (i.e., $m_{ij}(s) = 0, j \neq i$). That is to say, this decoupled system can only be composed of $p$ mutually independent controllable subsystems. Otherwise, assuming that a certain first-order branch (such as $x_k$) in the controllable and unobservable subsystem of this decoupled tree-like system is a common first-order branch in the controllable subsystems of $\bar{u}_i$ and $\bar{u}_{i'}$, and $path(\bar{u}_{i'}, y_i)$ can be constituted by the zero-order feedback branch of $x_k$ (denoted as $path_0(x_k, \bar{u}_i)$) after the closed-loop poles are assigned: $path(\bar{u}_{i'}, y_i) = path(\bar{u}_{i'}, x_k) + path_0(x_k, \bar{u}_i) + path(\bar{u}_i, y_i)$, so that the decoupled tree-like linear system is transformed into a non-decoupled system after the pole assignment. So there is the following theorem.

**Theorem 8.1.** For a regular equivalent decoupled tree-like controllable linear system (as shown in (29)), its closed-loop poles can be assigned while ensuring its decoupled state if and only if $m(s)$ in equation (156) is a diagonal matrix, in other words, this decoupled tree-like controllable system is only composed of $p$ mutually independent controllable subsystems.

**Theorem 8.2.** If a regular equivalent decoupled tree-like controllable system such as $\Sigma(\tilde{A}, \tilde{B}, \tilde{C})$ (as shown in (29)) is not composed of $p$ mutually independent controllable subsystems, its closed-loop poles cannot be assigned while ensuring its decoupled state, no matter what any coordinate transformation or regular static state feedback is used.

**Proof.** It is assumed that the system $\Sigma(\tilde{A}, \tilde{B}, \tilde{C})$ can be transformed into another regular equivalent decoupled tree-like controllable system $\Sigma(\tilde{A}', \tilde{B}', \tilde{C}')$ which is only composed of $p$ mutually independent controllable subsystems by a coordinate transformation $\bar{x} = T_1 \hat{x}$ and a static state feedback $(K_1, G_1)$ in turn. Therefore, there exists a regular static state feedback $(K_2, I)$ to assign the closed-loop poles of this decoupled system $\Sigma(\tilde{A}', \tilde{B}', \tilde{C}')$ while ensuring its decoupled state according to Theorem 8.1.

$$\bar{x} = T_1 \hat{x} \tag{222}$$
$$\upsilon = K_1 \bar{x} + G_1 \bar{\upsilon} \tag{223}$$
$$\bar{\upsilon} = K_2 \bar{x} + \tilde{\upsilon} \tag{224}$$
$$\upsilon = (K_1 + G_1 K_2) T_1 \hat{x} + G_1 \tilde{\upsilon} \tag{225}$$

The equation (225) shows that the closed-loop poles of $\Sigma(\tilde{A}, \tilde{B}, \tilde{C})$ can also be assigned by a regular static state feedback $((K_1 + G_1 K_2) T_1, G_1)$ while ensuring its decoupled state. Therefore, the system $\Sigma(\tilde{A}, \tilde{B}, \tilde{C})$ is only composed of $p$ mutually independent controllable subsystems according to Theorem 8.1. This is a contradiction, and so Theorem 8.2 holds.
□

Furthermore, Theorem 8.2 means that the regular equivalent decoupled tree-like controllable system $\Sigma(\tilde{A}, \tilde{B}, \tilde{C})$ is unique for a specific static state feedback $(F, G)$ to decouple the system $\Sigma(A, B, C)$. Moreover, if a regular equivalent decoupled tree-like controllable system $\Sigma(\tilde{A}, \tilde{B}, \tilde{C})$ is only composed of $p$ mutually independent controllable subsystems, then there exists one and only one independent static state feedback for each subsystem to assign its closed-loop poles.

## 9. Two examples for decoupling problem

**Example 9.1.** A tree-like controllable linear system is given by $\dot{x}_1 = u_1$; $\dot{x}_2 = u_2$; $\dot{x}_3 = x_4$; $\dot{x}_4 = x_5$; $\dot{x}_5 = x_6$; $\dot{x}_6 = u_3$; $\dot{x}_7 = x_8 + z, z = \alpha_1 x_2 + \alpha_2 x_3$; $\dot{x}_8 = x_9$; $\dot{x}_9 = u_4$, and $y_1 = x_1$; $y_2 = x_2$; $y_3 = x_1 + x_3$.

The signal flow graph corresponding to above system is drawn as shown in Fig.2. The decoupling framework is composed of three general paths such as $path(u_1, y_1)$, $path(u_2, y_2)$ and $path(u_3, y_3)$ whose general relative orders are $(1, 1, 4)$.

$$sy_1 = u_1; \quad sy_2 = u_2; \quad s^4 y_3 = u_3 + s^3 u_1$$

Therefore, control input $u_1$ needs to be compensated for one 3rd-order integrator-string, and only $path(u_4, x_7)$ is one 3rd-order integrator-string in the general decoupling framework. Then, $s^3 x_7 = u_4 + s^2 z$ means that input $z$ may need to be compensated for a certain chain of integrators.

If $\alpha_1 = 1, \alpha_2 = 0$, $s^3 x_7 = u_4 + su_2$, input $u_2$ needs to be compensated for one 1st-order integrator-string, but there is no independent integrator-string in the general decoupling framework. So, above system cannot be decoupled.

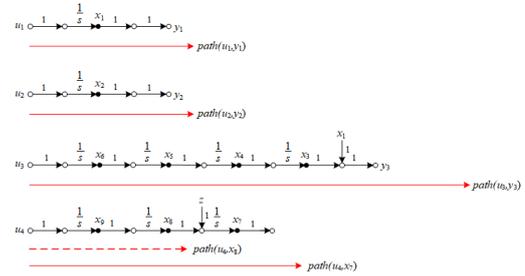

Fig.2. The signal flow graph and its general paths

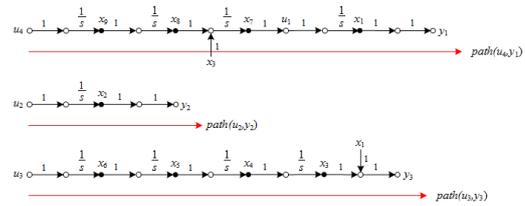

Fig.3. The decoupled signal flow graph

If $\alpha_1 = 0, \alpha_2 = 1$, $s^3 x_7 = u_4 + x_5$, there is one 3rd-order integrator-string $path(u_4, x_7)$ in the general decoupling framework. Let $u_1 = x_7$ (that is the 1st-type compensation) to decouple this system and its corresponding decoupled signal flow graph is shown in Fig.3.

$$s^4 y_1 = u_4 + x_5 = \upsilon_1; \quad sy_2 = u_2 = \upsilon_2;$$
$$s^4 y_3 = u_3 + u_4 + x_5 = \upsilon_3$$

Therefore, the static state feedback which decouple this system is shown as follows.

$$u_1 = x_7; \quad u_2 = \upsilon_2; \quad u_3 = \upsilon_3 - \upsilon_1; \quad u_4 = \upsilon_1 - x_5$$

The general relative orders of all general paths as shown in Fig.3 are $(4, 1, 4)$ and the relative orders of the decoupled system are the same. The decoupled system is a controllable and observable system and can be assigned the closed-loop poles while ensuring its decoupled state.

**Example 9.2.** Let $\Sigma(A, B, C)$ be given by $\dot{x}_1 = u_1$; $\dot{x}_2 = u_2$; $\dot{x}_3 = u_3$; $\dot{x}_4 = x_5$; $\dot{x}_5 = u_4$; $\dot{x}_6 = u_5$; $\dot{x}_7 = x_8$; $\dot{x}_8 = x_9$;

$\dot{x}_9 = x_{10}$ ; $\dot{x}_{10} = x_{11}$ ; $\dot{x}_{11} = u_6$ ; $\dot{x}_{12} = x_3 + x_{13}$ ; $\dot{x}_{13} = x_{14}$ ; $\dot{x}_{14} = x_1 + x_{15}$ ; $\dot{x}_{15} = x_{16}$ ; $\dot{x}_{16} = u_7$ ; $\dot{x}_{17} = x_{18}$ ; $\dot{x}_{18} = x_2 + x_{19}$ ; $\dot{x}_{19} = x_{20}$ ; $\dot{x}_{20} = u_8$ ; $\dot{x}_{21} = u_9$ ; $\dot{x}_{22} = u_{10}$ , and $y_1 = x_1$ ; $y_2 = x_2$ ; $y_3 = x_3$ ; $y_4 = x_3 - x_4$ ; $y_5 = x_6$ ; $y_6 = x_3 - x_4 + x_6 + x_7$ .

The signal flow graph corresponding to above system is drawn as shown in Fig.4. There are six general paths such as $path(u_1, y_1)$, $path(u_2, y_2)$, $path(u_3, y_3)$, $path(u_4, y_4)$, $path(u_5, y_5)$, and $path(u_6, y_6)$ whose general relative orders are $(1,1,1,2,1,5)$ in the decoupling framework.

$sy_1 = u_1$ ; $sy_2 = u_2$ ; $sy_3 = u_3$ ;
$s^2 y_4 = su_3 - u_4 = s(u_3 - s^{-1} u_4)$ ; $sy_5 = u_5$ ;
$s^5 y_6 = s^4 u_3 - s^3 u_4 + s^4 u_5 + u_6 = s^4(u_3 - s^{-1} u_4) + s^4 u_5 + u_6$

Then both mixed input $(u_3 - s^{-1} u_4)$ and control input $u_5$ need to be compensated for the 4th-order integrator-strings. $path(u_7, x_{13})$ and $path(u_8, x_{17})$ (or $path(u_7, x_{12})$ and $path(u_8, x_{17})$) can be selected as two general paths in the general decoupling framework.

If $path(u_7, x_{12})$ and $path(u_8, x_{17})$ are selected as two general paths, the following equations hold.
$s^5 x_{12} - u_7 = su_1 + s^3 u_3 = s(u_1 + s^2 u_3) = s^3(u_3 + s^{-2} u_1)$ ;
$s^4 x_{17} = su_2 + u_8$

It is obvious that two 1st-order and one 2nd-order (or one 1st-order and one 3rd-order) mutually independent integrator-strings are needed here, but only two mutually independent 1st-order integrator-strings are left in Fig.4. Therefore, the system $\sum(A, B, C)$ cannot be decoupled.

If $path(u_7, x_{13})$ and $path(u_8, x_{17})$ are selected as two general paths, the following equations hold.
$s^4 x_{13} = su_1 + u_7$ ; $s^4 x_{17} = su_2 + u_8$

So only two mutually independent 1st-order integrator-strings are needed. Let $u_1 = x_{21}$ and $u_2 = x_{22}$ (or $u_1 = x_{22}$, $u_2 = x_{21}$) to obtain two mutually independent 4th-order integrator-strings, and therefore, the system can be decoupled. Let $u_5 = x_{17}$ and $u_3 - s^{-1} u_4 = u_3 - x_5 = x_{13}$ (that is the 2nd-type compensation) to decouple this system and its decoupled signal flow graph is shown in Fig.5.

$s^2 y_1 = u_9 = \upsilon_1$ ; $s^2 y_2 = u_{10} = \upsilon_2$ ;
$y_3 = x_3$ ; $sy_3 = x_5 + x_{13}$ ; $s^5 y_3 = u_4 + x_{14} = \upsilon_3$ ;
$y_4 = x_3 - x_4$ ; $sy_4 = x_{13}$ ; $s^5 y_4 = u_7 + u_9 = \upsilon_4$ ;
$s^5 y_5 = u_8 + u_{10} = \upsilon_5$ ; $y_6 = x_3 - x_4 + x_6 + x_7$ ;
$sy_6 = x_{13} + x_{17} + x_8$ ; $s^5 y_6 = u_6 + u_7 + u_9 + u_8 + u_{10} = \upsilon_6$

Therefore, the static state feedback which decouple this system is shown as follows.
$u_1 = x_{21}$ ; $u_2 = x_{22}$ ; $u_3 = x_5 + x_{13}$ ; $u_5 = x_{17}$ ;
$u_4 = \upsilon_3 - x_{14}$ ; $u_6 = \upsilon_6 - \upsilon_4 - \upsilon_5$ ; $u_7 = \upsilon_4 - \upsilon_1$ ;
$u_8 = \upsilon_5 - \upsilon_2$ ; $u_9 = \upsilon_1$ ; $u_{10} = \upsilon_2$

The general relative orders of all general paths as shown in Fig.5 are $(2,2,5,2,5,5)$ and the relative orders of the decoupled system are $(2,2,2,5,5,5)$.

Only 1st-order branch $x_{12}$ (whose input node is connected to the outputs of two 1st-order branches such as $x_3$ and $x_{13}$, and is further connected to $y_3$, $y_3^{(1)}$, $y_4$, $y_4^{(1)}$, $y_6$ and $y_6^{(1)}$) is controllable and unobservable in the regular equivalent decoupled tree-like controllable system, and so it is a common 1st-order branch in the controllable subsystems of $\upsilon_3$, $\upsilon_4$ and $\upsilon_6$. Therefore, the regular equivalent decoupled tree-like controllable system corresponding to the system $\sum(A, B, C)$ is not composed of six mutually independent controllable subsystems and cannot be assigned the closed-loop poles while ensuring its decoupled state.

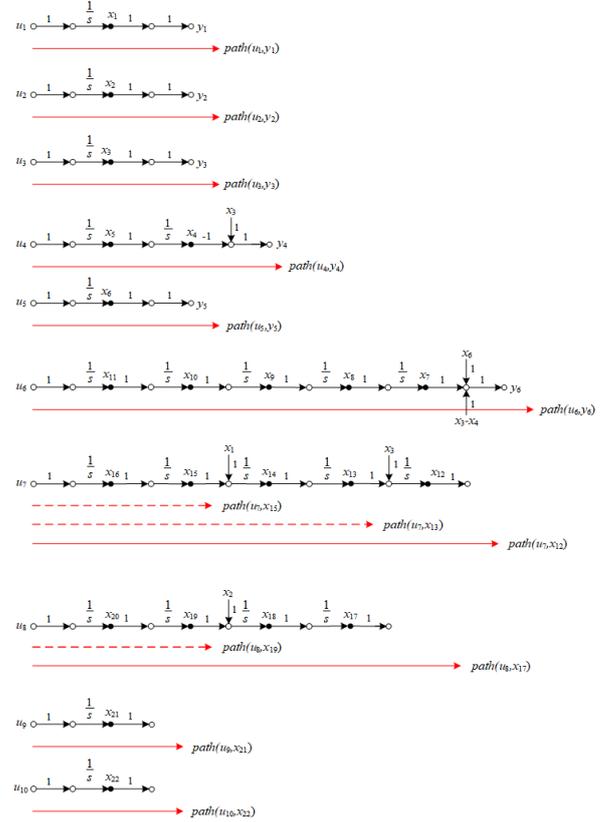

Fig.4. The signal flow graph and its general paths

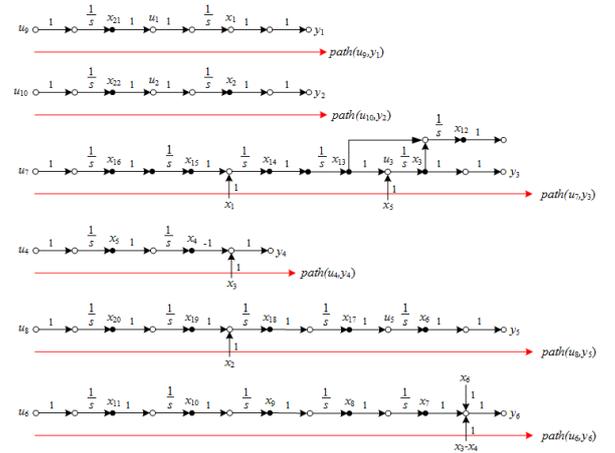

Fig.5. The decoupled signal flow graph

If let $u_5 = x_{13}$ and $u_3 - x_5 = x_{17}$ to decoupled this system, only 1st-order branch $x_{12}$ (whose input node is connected to the outputs of two 1st-order branches such as $x_3$ and $x_{13}$, and is further connected to $y_3$, $y_4$, $y_5^{(1)}$, $y_6$ and $y_6^{(1)}$) is also controllable and unobservable in another regular equivalent decoupled tree-like controllable system, and so it is a common 1st-order branch in the controllable subsystems of $\upsilon_3$, $\upsilon_4$, $\upsilon_5$ and $\upsilon_6$. Therefore, the conclusion about the pole assignment for this regular equivalent decoupled tree-like controllable system is the same as above.

## 10. Conclusions

This paper has established the necessary and sufficient condition for the solvability of a nonregular Morgan's problem by a nonregular static state feedback and has also obtained the static state feedback expression for decoupling system. Moreover, the stability problem of a decoupled system has also solved. The internal stability of the uncontrollable subsystem of a decoupled system is ensured by assigning its closed-loop poles during the self-compensation process of all auxiliary inputs, but the stability of the controllable and unobservable subsystem of a decoupled system depends on whether the controllable subsystem of the decoupled system is composed of $p$ mutually independent controllable subsystems.

Some interesting potential areas of future research arise from the results obtained in this paper. In particular, the question of extending the theory to the time-varying linear or non-linear situation is of considerable interest.

**Acknowledgements**

The author would like to thank Professor Wen Dingdou and Dr. Zhang Yang for their helpful discussion while preparing this work, and would like to especially thank Dr. Zhang Yang for drawing all figures in this paper.